\documentclass[12pt]{article}
\usepackage{latexsym,amssymb,
amsthm,amsfonts,amsmath, epsfig,psfrag,array}

\begin{document}

\newtheorem{lemma}{Lemma}
\newtheorem{cor}{Corollary}
\newtheorem{theorem}{Theorem}
\renewcommand{\proofname}{Proof}
\newtheorem{property}{Property}
\newtheorem{prop}{Proposition}[subsection]

\def \H{{\mathbb H}}
\def \R{{\mathbb R}}
\def \E{{\mathbb E}}
\def \S{{\mathbb S}}
\def \G{{\mathcal G}}

\def \l{\langle }
\def \r{\rangle }
\def \[{[ }
\def \]{] }
\def \d{D\,}
\def \sign{\text{\,sign\,}}
\def \conv{\text{\,conv\,}}
\def \o{\overline}
\def \wt{\widetilde}
\def \wh{\widehat}
\def \a{\alpha}
\def \b{\beta}
\def \dim{\text{\,dim}}

\begin{center}
{
\LARGE  \bf
 On hyperbolic Coxeter polytopes \\
\smallskip
with mutually intersecting facets\\
}

\vspace{19pt}
{\large Anna~Felikson\footnote{Partially supported by grant of President of Russia NS-5666.2006.1 }}
\qquad\qquad
{\large Pavel~Tumarkin\footnote{Partially supported by grants of President of Russia MK-6290.2006.1 and NS-5666.2006.1}}

\smallskip
felikson@mccme.ru\qquad\qquad\quad
pasha@mccme.ru
\vspace{15pt}

\begin{tabular}{c}
{\small Independent University of Moscow}\\
{\small B. Vlassievskii 11, 119002 Moscow, Russia}
\end{tabular}

\end{center}
\vspace{12pt}

\noindent
\begin{center}
\raisebox{20pt}{}{\parbox{12cm}{{\bf Abstract.\ }\small
We prove that, apart from some well-known low-dimensional examples,
any compact hyperbolic Coxeter polytope has a pair of disjoint
facets. This is one of very few known general results concerning
combinatorics of compact hyperbolic Coxeter polytopes.
We also obtain a similar result for simple non-compact polytopes.
} }

\end{center}

\bigskip

%\tableofcontents

%\vspace{2pt}

\section{Introduction.}

A Coxeter polytope in the spherical, hyperbolic or
Euclidean space is a polytope whose dihedral angles are all
integer submultiples of $\pi$.
These polytopes are very important among
acute-angled polytopes since
a group generated by reflections with respect to the facets of a Coxeter
polytope is discrete.
On the other hand, a fundamental chamber of any (finitely generated)
discrete reflection group in these spaces is a Coxeter polytope.

Already in 1934, H.~S.~M.~Coxeter~\cite{C} proved that any spherical
Coxeter polytope (containing no pair of opposite points of the sphere)
is a simplex and any compact Euclidean Coxeter polytope is either a
simplex or a direct product of simplices.

However,  hyperbolic Coxeter polytopes are still far
from being classified.
It was proved by E.~Vinberg~\cite{V1} that no compact hyperbolic
Coxeter polytope exists in dimensions $d\ge 30$; M.~Prokhorov~\cite{996}
and A.~Khovanskij~\cite{Kh} proved that no hyperbolic
Coxeter polytope of finite volume exists in dimensions $d\ge 996$.
These bounds do not look sharp: the examples are known only up to
dimension $8$ in compact case and up to dimension $21$ in the
non-compact case.

\begin{center}
\begin{figure}[!ht]
\begin{center}
\psfrag{8}{\small $8$}
\psfrag{10}{\small $10$}
\epsfig{file=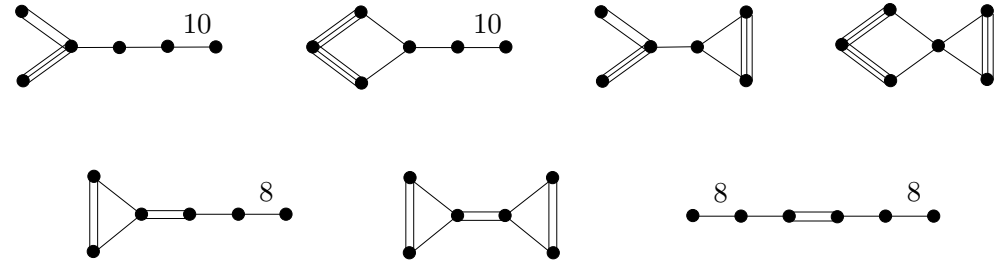,width=0.95\linewidth}
\caption{Coxeter diagrams of Esselmann polytopes (or
{\it Esselmann diagrams} for short).}
\label{Ess}
\end{center}
\end{figure}
\end{center}

\begin{center}
\begin{figure}[!hb]
\begin{center}
\epsfig{file=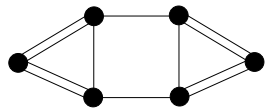,width=0.25\linewidth}
\caption{Coxeter diagram of the unique simple non-compact
Coxeter $d$-polytope that has $d+2$ facets and has no
pair of disjoint facets.}
\label{t}
\end{center}
\end{figure}
\end{center}

Besides the restriction on the dimension and some series of
examples,
%~\cite{M},~\cite{ImH},~\cite{Al},
there exists a classification of hyperbolic Coxeter
polytopes of certain combinatorial types.
More precisely, compact simplices were classified by
F.~Lan\-n\'er~\cite{L}, and non-compact simplices were classified
by several authors (see e.g.~\cite{Bou} or~\cite{V67}).
Simplicial prisms were listed by I.~Kaplinskaja~\cite{K};
F.~Esselmann~\cite{Ess} obtained the classification of the
remaining compact hyperbolic Coxeter $d$-dimensional polytopes
with $d+2$ facets. These consist of seven $4$-dimensional
polytopes with mutually intersecting facets (we call these
polytopes {\it
Esselmann polytopes} and reproduce the list in
Fig.~\ref{Ess}).
P.~Tumarkin~\cite{T} classified those non-compact hyperbolic Coxeter
$d$-dimensional polytopes with $d+2$ facets that do not have disjoint facets.
The only simple polytope from this list is shown in Figure~\ref{t}
and is of dimension $4$.

\bigskip

This paper is devoted to the proof of the following theorem:

%\begin{theorem}
%\label{nodots}

\bigskip
\noindent
{\bf Theorem~A.}
{\it
Let $P$ be a compact hyperbolic Coxeter $d$-dimensional polytope.
If $d>4$ then $P$ has a pair of disjoint facets.

If $d\le 4$ and $P$ has no pair of disjoint facets
then $P$ is either a simplex or one of the seven Esselmann
polytopes. }
%\end{theorem}

\bigskip
\noindent
A $d$-dimensional polytope is {\it simple} if any vertex of $P$ is contained in
exactly $d$ facets, or equivalently, facets of $P$ at each vertex
are in general position. The classification of spherical polytopes
implies that any compact hyperbolic Coxeter polytope is simple.
While proving
Theorem~A, we slightly change the proof to obtain a similar result
concerning simple non-compact hyperbolic Coxeter polytopes of
finite volume.

\bigskip
\noindent
{\bf Theorem~B.}
{\it
Let $P$ be a simple non-compact hyperbolic Coxeter $d$-dimen\-sional
polytope of finite volume.
If $d>9$ then $P$ has a pair of disjoint facets.

If $d\le 9$ and $P$ has no pair of disjoint facets
then $P$ is either a simplex or the $4$-dimensional polytope shown
in Fig.~\ref{t}. }

\bigskip

\noindent
The paper is organized as follows.
In Section~\ref{prelim} we recall some information about Coxeter polytopes.
In Section~\ref{technical tools} we introduce some technical tools
we use for proving the theorems.
Section~\ref{proof} is devoted to the proof of the theorems.

We prove both theorems A and B simultaneously.
The proof is by induction on the dimension $d$.
The most general case is $d\ge 9$.
In this case the proof is by examination of the combinatorics of
the Coxeter diagram $\Sigma(P)$ while making use of a recent
result of
D.~Allcock (Theorem~\ref{All}).
Some minor technical refinements generalize the proof to
$d\ge 7$ (see Section~\ref{d9}).

The small dimensions are considered in
Sections~\ref{d2,3},~\ref{d4},~\ref{d5} and~\ref{d6}.
In dimensions $d=2$ and $3$ the argument is purely combinatorial
(Lemma~\ref{2,3}). In dimensions from $4$ to $6$ the proof also
 uses a computational
technique developed in Section~\ref{technical tools}
based on the notion of local determinants.

%Note also, that any 1-dimensional polytope is a segment,
%so it always has a pair of disjoint facets (i.e. two distinct vertices).

The authors would like to thank the referee for many helpful comments
and suggestions.  We are grateful to D.~Panov for pointing out to us the
result of D.~Allcock, and to D.~Allcock who
introduced us to his paper~\cite{Al}, which was not published yet.

\section{Preliminaries.}
\label{prelim}
In this section we list the essential facts about Coxeter diagrams and
Gale diagrams.
Concerning Coxeter diagrams we follow mainly~\cite{V2} and~\cite{29}.
For details about Gale diagrams see~\cite{G}.
At the end of the section we recall results of R.~Borcherds~\cite{B} and
D.~Allcock~\cite{Al} concerning Coxeter faces of Coxeter polytopes.

In what follows we write {\it $d$-polytope} instead of
``$d$-dimensional polytope'', {\it $k$-face} instead of ``$k$-dimensional
face'' and {\it facet} instead of ``face of codimension one''.

\subsection{Coxeter diagrams.}

{\bf 1.}
An abstract Coxeter diagram $\Sigma$ is a finite $1$-dimensional
simplicial complex with weighted edges, where weights $w_{ij}$ are
positive, and if $w_{ij}<1$ then $w_{ij}=\cos \frac{\pi}{k}$ for
some integer $k\ge 3$.
A {\it sub\-dia\-gram} of $\Sigma$ is a subcomplex with the same weights
as in $\Sigma$.
The {\it order}  $|\Sigma|$ is the number of vertices
of the diagram $\Sigma$.

If $\Sigma_1$ and $\Sigma_2$ are sub\-dia\-grams of an abstract
Coxeter diagram $\Sigma$,
we denote by $\l \Sigma_1,\Sigma_2\r$ a sub\-dia\-gram of $\Sigma$
spanned by the vertices of $\Sigma_1$ and $\Sigma_2$.

Given an abstract Coxeter diagram $\Sigma$ with vertices
$v_1,\dots,v_n $ and weights $w_{ij}$, we construct a symmetric
$n\times n$ matrix $M(\Sigma)=(c_{ij})$, where $c_{ii}=1$,
$c_{ij}= -w_{ij}$ if  $v_i$ and $v_j$ are  adjacent, and
$c_{ij}=0$ otherwise.
By determinant, rank and signature of $\Sigma$ we mean the
determinant, the rank and the signature of $M(\Sigma)$.

We can suppress the weights but indicate the same information by
labelling the edges of a Coxeter diagram in the following way: if
the weight $w_{ij}$ equals $\cos\frac{\pi}{m}$, $v_i$ and $v_j$
are joined by an $(m-2)$-fold edge or a simple edge labelled by
$m$; if $w_{ij}=1$, $v_i$ and $v_j$ are joined by a bold edge; if
$w_{ij}>1$, $v_i$ and $v_j$ are joined by a dotted edge labelled
by $w_{ij}$ (or without any label).

We write $\[v_i,v_j\]=m$ if $w_{ij}$ equals
$\cos\frac{\pi}{m}$, and $\[v_i,v_j\]=\infty$ if $w_{ij}\ge 1$.

By a {\it multiple} edge we mean an edge of weight
$\cos\frac{\pi}{m}$ for $m\ge 4$. By a {\it multi-multiple}
edge we mean an edge of weight $\cos\frac{\pi}{m}$ for $m\ge 6$.

\vspace{7pt}

%\begin{defin}
%$\bullet$

An abstract Coxeter diagram $\Sigma$ is {\it elliptic} if
$M(\Sigma)$ is positive definite. An $n\times n$ matrix is called
{\it indecomposable} if it cannot be transformed to a
block-diagonal one by simultaneous permutations of columns and
rows. Clearly, connected components of $\Sigma$ correspond to
indecomposable components of $M(\Sigma)$. A diagram $\Sigma$ is
{\it parabolic} if any indecomposable component of $M(\Sigma)$ is
degenerate and positive semidefinite; a connected diagram $\Sigma$
is a {\it Lann\'er} diagram if $\Sigma$ is indefinite but any
proper sub\-dia\-gram of $\Sigma$ is elliptic; a connected diagram
$\Sigma$ is a {\it quasi-Lann\'er} diagram if $\Sigma$ is not a
Lann\'er diagram, $\Sigma$ is indefinite, but any proper
sub\-dia\-gram of $\Sigma$ is either elliptic or parabolic;
$\Sigma$  is {\it superhyperbolic} if its negative inertia index
is greater than $1$.

\vspace{7pt}

The list of connected elliptic and parabolic diagrams with
their standard notation is contained in~\cite[Tables~1,2]{29}.
See also~\cite[Tables~3,4]{29} for the lists of Lann\'er and
quasi-Lann\'er diagrams.
We need the following properties of these lists:

\begin{itemize}
%\vspace{4pt}
%\noindent
\item[$\bullet\,$] %\quad
there are finitely many Lann\'er diagrams of order greater than
$3$, and the maximal order of a Lann\'er diagram is $5$;
\item[$\bullet\;$]  %\quad
any Lann\'er diagram of order $5$ contains a
sub\-dia\-gram of the type $H_4$ or~$F_4$;
\item[$\bullet\;$]  %\quad
any Lann\'er diagram of order $4$ contains a
sub\-dia\-gram of the type $H_3$ or~$B_3$;
\item[$\bullet\;$]  %\quad
any Lann\'er diagram of order $3$ contains a
multiple edge;
\item[$\bullet\;$]  %\quad
 Lann\'er diagrams of order greater than $3$ contain no
multi-multiple edges;
\item[$\bullet\;$]
any quasi-Lann\'er diagram of order $n$ contains a connected
parabolic sub\-dia\-gram of order $n-1$.

\end{itemize}

\bigskip

\noindent
{\bf 2.}
It is convenient to describe Coxeter polytopes by their Coxeter diagrams.
Let $P$ be a Coxeter polytope with facets $f_1,\dots,f_r$.
The Coxeter diagram $\Sigma(P)$ of the polytope $P$
is a diagram with vertices $v_1,\dots,v_r$; two edges $v_i$ and
$v_j$ are not joined if the hyperplanes spanned by $f_i$ and $f_j$
are orthogonal;
 $v_i$ and $v_j$ are joined by an edge with weight
$$w_{ij}=
\begin{cases}
\cos \frac{\pi}{k},&\text{ if  $f_i$ and $f_j$ form a  dihedral angle
$\frac{\pi}{k}$;}\\
1,&\text{ if $f_i$ is parallel to $f_j$;}  \\
\cosh \rho,&\text{ if $f_i$ and $f_j$ diverge and $\rho$ is the distance
from $f_i$ to $f_j$.}
\end{cases}
$$

\noindent
If $\Sigma=\Sigma(P)$, then $M(\Sigma)$ coincides with the Gram matrix
of outer unit normals to the facets of $P$ (referring to the
standard model of hyperbolic $d$-space in $\R^{d,1}$).

It is shown in~\cite{V2} that a Coxeter diagram $\Sigma(P)$
of a compact $d$-dimensi\-onal hyperbolic polytope $P$ is a
connected diagram of signature $(d,1)$ without parabolic
sub\-dia\-grams.
In particular, $\Sigma(P)$ contains no bold edge, and any indefinite
sub\-dia\-gram contains a Lann\'er diagram. Moreover, it is shown
there that any compact hyperbolic Coxeter $d$-polytope $P$ is
simple, and elliptic sub\-dia\-grams of $\Sigma(P)$ are in
one-to-one correspondence with faces of $P$: a $k$-face $F$
corresponds to an elliptic sub\-dia\-gram of order $d-k$ whose
vertices correspond to the facets of $P$ containing $F$.

It is also shown in~\cite{V2} that if $\Sigma(P)$ is a Coxeter diagram
of a non-compact hyperbolic $d$-polytope $P$, then for any ideal
vertex $V$ of $P$ (i.e. $V$ lies at the boundary of the hyperbolic
space) the vertices of
$\Sigma(P)$ corresponding to facets containing $V$ compose a
parabolic diagram of rank
$d-1$, and any parabolic sub\-dia\-gram of $\Sigma(P)$ may be enlarged
to some parabolic sub\-dia\-gram of rank $d-1$. In particular, if
$P$ is simple then any parabolic sub\-dia\-gram $S$ of $\Sigma(P)$ is
connected and has rank $d-1$, i.e. $S$ has order $d$. Clearly, any
indefinite sub\-dia\-gram of $\Sigma(P)$ contains either a Lann\'er or
quasi-Lann\'er diagram.

As an easy corollary, we have the following statement.

\begin{prop}
\label{small}
Let $P$ be a simple hyperbolic Coxeter $d$-polytope.
Then $\Sigma(P)$ contains either a Lann\'er or quasi-Lann\'er
diagram, and $\Sigma(P)$ does not contain parabolic diagrams of
order less than $d$.

\end{prop}

\begin{lemma}
\label{same}
Let $\Sigma(P)$ be a Coxeter diagram of a hyperbolic Coxeter $d$-polytope
$P$ of finite volume. Then no proper sub\-dia\-gram of $\Sigma(P)$ is a diagram
of a hyperbolic Coxeter $d$-polytope of finite volume.

\end{lemma}

\begin{proof}
Suppose that a proper sub\-dia\-gram $\Sigma \subset \Sigma(P)$
is a diagram of a Coxeter $d$-polytope of finite volume.
The vertices of $\Sigma$ determine a polytope $P'$.
Denote by $G_P$ and $G_{P'}$ the groups generated by
reflections with respect to the facets of $P$ and $P'$
respectively.
The group $G_{P'}$ is a
subgroup of $G_P$. Since $P'$ is of finite volume,  $G_{P'}$ has a
finite index in $G_P$. At the same time, the number of facets of
$P$ is more than $P'$ has. This contradicts the main result
of~\cite{FT} which claims that if $P$ and $P'$ are finite volume
Coxeter polytopes in $\H^n$ or $\E^n$, $G_{P}$ and $G_{P'}$ are the
groups generated by reflections in the facets of $P$ and $P'$
respectively, and $G_{P'}\subseteq G_P$ is a finite index
subgroup, then the number of facets of $P$ does not exceed the
number of facets of $P'$.

\end{proof}

\begin{cor}
\label{ql}
If a Coxeter diagram of a simple Coxeter polytope $P$ contains
a quasi-Lann\'er sub\-dia\-gram then $P$ is a simplex.

\end{cor}

\begin{proof}
Any quasi-Lann\'er diagram of order $d+1$ is a Coxeter diagram of non-compact
hyperbolic Coxeter $d$-dimensional simplex of finite volume (see
e.g.~\cite{V2}).
Suppose that $P$ is not a simplex.
Lemma~\ref{same} implies that if $P$ is a $d$-polytope of finite volume then
$\Sigma(P)$ contains no quasi-Lann\'er sub\-dia\-grams of order $d+1$.
Clearly, $\Sigma(P)$ does not contain any quasi-Lann\'er sub\-dia\-gram
of order greater than $d+1$.
Further, since $P$ is simple, any connected parabolic sub\-dia\-gram of
$\Sigma(P)$ should have order $d$, so $\Sigma(P)$ contains no quasi-Lann\'er
sub\-dia\-gram of order less than $d+1$, either.

\end{proof}

\bigskip

\subsection{Gale diagrams and missing faces.}
We do not use the content of this section throughout the paper
except for the proof of the Theorems~A and~B for $4$-polytopes.

Every combinatorial type of simple $d$-polytope with
$d+k$ facets can be represented by its {\it Gale diagram} $\G$.
This consists of $d+k$ points $a_1,\dots,a_{d+k}$ on
$(k-2)$-dimensional unit sphere $\S^{k-2}\subset \R^{k-1}$
centered at the origin. Each point $a_i$ corresponds to a facet
$f_i$ of $P$. The combinatorial type of a convex polytope can be
read off from the Gale diagram in the following way: for any
$J\subset \{1,\dots,d+k\}$ the intersection of facets $\{f_j \,|\,
j\in J \}$ is a proper (that is, non-empty) face of $P$ if and
only if the origin is contained in the interior of $\conv\{a_{j}
\,|\, j\notin J\}$ (where $\conv X$ is a convex hull of the set
$X$).

The points $a_1,\dots,a_{d+k}\in \S^{k-2}$ compose a Gale diagram of
some $d$-dimensional polytope $P$ with $d+k$ facets if and only if
every open half-space $H^+$ in $\R^{k-1}$ bounded by a hyperplane
$H$ through the origin contains at least two of the points
$a_1,\dots,a_{d+k}$.

Notice that the definition of Gale diagram introduced above concerns
simple polytopes only, and it is ``dual" to the standard one (see,
for example,~\cite{G}): usually Gale diagram is defined in terms
of vertices of polytope instead of facets. Notice also that the
definition above takes simplices out of consideration: usually one
means the origin of $\R^1$ with multiplicity $d+1$ by the Gale
diagram of a $d$-simplex, however we exclude the origin since we
consider simple polytopes only, and the origin is not contained in
$\G$ for any simple polytope except the simplex.

\medskip

Let $P$ be a simple polytope.
The facets $f_1,\dots,f_m$ of $P$ compose a {\it missing face} of $P$
if $\bigcap\limits_{i=1}^m f_i=\emptyset$ but any proper subset of
$\{f_1,\dots,f_m\}$ has a non-empty intersection.

\begin{lemma}
\label{missing}
Let $P$ be a simple $d$-polytope with $d+k$ facets $\{f_i\}$, let
$\G=\{a_i\}\subset\S^{k-2}$ be a Gale diagram of $P$, and let
$I\subseteq\{1,\dots,d+k\}$.
Then the set $M_I=\{f_i\,|\,i\in I\}$ is a missing face of
$P$ if and only if the following two conditions hold:
\begin{itemize}
\item[(1)] there exists a hyperplane $H$ through the origin separating
the set $\wh M_I=\{a_i\,|\,i\in I\}$ from the rest points of $\G$;
\item[(2)] for any proper subset $J\subset I$ no hyperplane through
the origin separates the set $\wh M_J=\{a_i\,|\,i\in J\}$ from the
remaining points of $\G$.

\end{itemize}
\end{lemma}

\begin{proof}
Suppose first that both conditions hold.
Since $P$ is simple, (1) implies that $\conv(G\setminus \wh M_I)$
does not contain the origin, so
$\bigcap\limits_{i\in I} f_i=\emptyset$.
If $\bigcap\limits_{i\in J} f_i$ is also empty for some
$J\subsetneq I$, we obtain that $\conv(\G\setminus\wh M_J)$ does
not contain the origin, so there exists a hyperplane $H$ through
the origin such that $\G\setminus\wh M_J$ is contained in one of
halfspaces $H^+$ and $H^-$, say $H^+$. Then $\G\cap H^-$ is a
subset of $\wh M_J$, i.e. some subset of $\wh M_J$ is separated by
a hyperplane through the origin, which contradicts (2).

Now suppose that $M_I$ is a missing face. Then there exists a
hyperplane $H$ through the origin such that $\G\setminus\wh M_I$
is contained in a halfspace $H^+$. Since $P$ is simple, we may
assume that $\G\cap H=\emptyset$. To prove (1) suppose the
contrary, i.e. $a_{i_0}\in H^+$ for some $i_0\in I$.
Then $\G\setminus\wh M_{I\setminus i_0}$ is also
contained in $H^+$, that means that
$\bigcap\limits_{i\in {I\setminus i_0}}\!\!\!f_i$ is empty in
contradiction to the definition of missing face.
To prove (2) notice that if some hyperplane $H_J$ separates
$\wh M_J$ for some $J\subsetneq I$ then $\bigcap\limits_{i\in J}
f_i=\emptyset$, which also contradicts the definition of missing
face.

\end{proof}

\bigskip

Suppose that $P$ is a simple hyperbolic Coxeter polytope.
The definition of missing face implies that for any
Lann\'er or quasi-Lann\'er sub\-dia\-gram $L\subset\Sigma(P)$
the facets corresponding to $L$ compose a missing face of $P$ (and
any missing face of $P$ corresponds to some Lann\'er or
quasi-Lann\'er diagram in $\Sigma(P)$).

\subsection{Faces of Coxeter polytopes.}

Let $P$ be a hyperbolic Coxeter $d$-polytope, and denote by $\Sigma(P)$
its Coxeter diagram. Let $S_0$ be an elliptic sub\-dia\-gram of
$\Sigma(P)$.
By~\cite[Th.~3.1]{V2}, $S_0$ corresponds to a face of $P$ of dimension
$d-|S_0|$. Denote this face by $P(S_0)$. $P(S_0)$ itself is an
acute-angled polytope, but it might not be a Coxeter polytope.
R.~Borcherds proved the following
sufficient condition for $P(S_0)$ to be a Coxeter polytope.

\begin{theorem}[\cite{B}, Ex. 5.6]
\label{bor}
Suppose $P$ is a Coxeter polytope with diagram $\Sigma(P)$,
and $S_0\subset \Sigma(P)$ is an elliptic sub\-dia\-gram that has no
$A_n$ or $D_5$ component. Then $P(S_0)$ itself is a Coxeter polytope.

\end{theorem}

Facets of $P(S_0)$ correspond to those vertices
that together with $S_0$ comprise an elliptic or positive
semidefinite subdiagram of $\Sigma(P)$.
The following result of D.~Allcock shows how to compute dihedral
angles of $P(S_0)$.

Let $a$ and $b$ be the facets of $P(S_0)$ coming from facets $A$ and $B$ of
$P$, i.e. $a=A\cap P(S_0)$ and $b=B\cap P(S_0)$. Denote by $v_A$
and $v_B$ the nodes of $\Sigma(P)$ corresponding to the facets $A$
and $B$. We say that a node of $\Sigma(P)$ {\it attaches} to $S_0$
if it is joined with some nodes of $S_0$ by edges of any type.
Then the angles of $P(S_0)$ can be computed in the following way.

\begin{theorem}[\cite{Al}, Th. 2.2]
\label{All}
Under the hypotheses of Theorem~\ref{bor},
\begin{itemize}
\item[(1)]
If neither $v_A$ nor $v_B$ attaches to $S_0$, then $\angle ab=\angle AB$.
\item[(2)]
If just one of $v_A$ and $v_B$ attaches to $S_0$, say to the component
$S_0^i$, then
\begin{itemize}
\item[(a)]
if $A\perp B$ then $a\perp b$;
\item[(b)]
if $v_A$ and $v_B$ are joined by a simple edge, and adjoining
$v_A$ and $v_B$ to $S_0^i$ yields a diagram $B_k$ (resp. $D_k$,
$E_8$ or $H_4$) then $\angle ab=\pi/4$ (resp. $\pi/4$, $\pi/6$ or
$\pi/10$);
\item[(c)]
otherwise, $a$ and $b$ do not meet.
\end{itemize}
\item[(3)]
If $v_A$ and $v_B$ attach to different components of $S_0$, then
\begin{itemize}
\item[(a)]
if $A\perp B$ then $a\perp b$;
\item[(b)]
otherwise, $a$ and $b$ do not meet.
\end{itemize}
\item[(4)]
If $v_A$ and $v_B$ attach to the same component of $S_0$, say $S_0^i$,
then
\begin{itemize}
\item[(a)]
if $A$ and $B$ are not joined and $S_0^i\cup \{A,B\}$ is a diagram
$E_6$ (resp. $E_8$ or $F_4$) then $\angle ab=\pi/3$ (resp. $\pi/4$
or $\pi/4$);
\item[(b)]
otherwise, $a$ and $b$ do not meet.
\end{itemize}
\end{itemize}

\end{theorem}

Let $w\in \Sigma(P)$ be a {\it neighbor} of $S_0$,
so that $w$ attaches to $S_0$ by some edges.
We call $w$ a {\it good neighbor} if $\l S_0, w\r$ is either an
elliptic diagram or a positive semidefinite diagram, and
{\it bad} otherwise.
We denote by $\overline S_0$ the sub\-dia\-gram of $\Sigma(P)$
consisting of vertices corresponding to facets of $P(S_0)$. The
diagram $\overline S_0$ is spanned by good neighbors of $S_0$ and
by all vertices not joined to $S_0$ (in other words, $\overline
S_0$ is spanned by all vertices of $\Sigma(P)\setminus S_0$ except
bad neighbors of $S_0$). If $P(S_0)$ is a Coxeter polytope, denote
its Coxeter diagram by $\Sigma_{S_0}$.

\begin{cor}
\label{cor_All}
Suppose that $P(S_0)$ is a Coxeter polytope.
\begin{itemize}
\item[(a)]
If $S_0$ has no good neighbors then $\overline S_0=\Sigma_{S_0}$.
In particular, this always holds for $S_0=H_4$ and $G_2^{(m)}$
where $m\ge 7$, for $S_0=H_4$ if $d>4$, and for $S_0=G_2^{(6)}$ if
$d>3$.

\item[(b)]
If $S_0=B_n$, $n\ge 2$, and $\Sigma_{S_0}$ contains a sub\-dia\-gram $S$ of
the type $H_4$ or $F_4$, then $S$ is contained in $\overline S_0$, too.

\end{itemize}
\end{cor}

\begin{proof}
To prove (a) one should only notice that all neighbors of diagrams
listed in item (a) (except for $F_4$ and $G_2^{(6)}$)  are bad.
Any good neighbor of $F_4$ or $G_2^{(6)}$ leads to a parabolic sub\-dia\-gram of
$\Sigma(P)$ of order $4$ and $3$ respectively, which contradicts
Prop.~\ref{small} in case of $d>4$ and $d>3$.

Item (b) follows immediately from Theorem~\ref{All}.

\end{proof}

Notice also that any face of a simple polytope is a simple
polytope itself. In particular, if $P$ is simple then for any
elliptic sub\-dia\-gram $S\subset\Sigma$ the polytope $P(S)$ is also
simple.

\section{Technical tools.}
\label{technical tools}

\subsection{Local determinants.}

Let $\Sigma$ be a Coxeter diagram, and let $T$ be a sub\-dia\-gram
of $\Sigma$ such that $\det(\Sigma\setminus T)\ne 0$.
A {\it local determinant} of $\Sigma$ on a sub\-dia\-gram $T$ is
$$\det(\Sigma,T)=\frac{\det \Sigma}{\det(\Sigma\!\setminus\! T)}.$$

\begin{prop}[\cite{V1}, Prop.~12]
\label{loc_sum}
If a Coxeter diagram $\Sigma$ consists of two sub\-dia\-grams
$\Sigma_1$ and $\Sigma_2$ having a unique vertex $v$ in common,
and no vertex of $\Sigma_1\setminus v$ attaches to
$\Sigma_2\setminus v$, then
$$ \det(\Sigma,v)=\det(\Sigma_1,v)+\det(\Sigma_2,v)-1.$$

\end{prop}

\begin{prop}[\cite{V1}, Prop.~13 ]
\label{loc_product}
If a Coxeter diagram $\Sigma$ is spanned by two disjoint sub\-dia\-grams
$\Sigma_1$ and $\Sigma_2$ joined by a unique edge $v_1v_2$ such that $[v_1,v_2]=m$,
then
$$\det(\Sigma,\l v_1,v_2\r )=\det(\Sigma_1,v_1) \det(\Sigma_2,v_2) - \cos^2\frac{\pi}{m}.$$

\end{prop}

\begin{figure}[!h]
\begin{center}
\psfrag{a}{ $a$}
\psfrag{b}{ $b$}
\psfrag{c}{ $c$}
\psfrag{v}{ $v$}
\epsfig{file=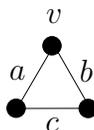,width=0.1\linewidth}
\caption{Diagram $L_{a,b,c}$}
\label{d_abc}
\end{center}
\end{figure}

Denote by $L_{abc}$ a Lann\'er diagram of order $3$
containing sub\-dia\-grams of the dihedral groups $G_2^{(a)}$,
$G_2^{(b)}$ and $G_2^{(c)}$.
Let $v$ be the vertex of $L_{a,b,c}$ that does not belong to
$G_2^{(c)}$.
Denote by $\d(a,b,c)$ the local determinant $\det(L_{a,b,c},v)$,
see Fig.~\ref{d_abc}.

It is easy to check (see e.g.~\cite{V1}) that
$$
\d(a,b,c)=
1-\frac{\cos^2(\pi/a)+\cos^2(\pi/b)+2\cos(\pi/a)\cos(\pi/b)\cos(\pi/c)}
{\sin^2(\pi/c)}.
$$

Notice that $|\d(a,b,c)|$ is an increasing function on each
of $a,b,c$ tending to
infinity while $c$ tends to infinity.

\subsection{Lists $L_\a(S_0,d)$, $L_\b(S_0,d)$ and $L'(\Sigma,C,d)$.}

\begin{lemma}
\label{pr1}
Let $P$ be a simple Coxeter $d$-polytope with mutually intersecting
facets, and assume that $P$ is not a simplex. Let $S_0$ be a
connected elliptic sub\-dia\-gram of
$\Sigma(P)$ such that \\
{\sc (i)} $|S_0|< d$ and $S_0\ne A_n$, $D_5$.\\
{\sc (ii)} $S_0$ has no good neighbors in $\Sigma(P)$.\\
{\sc (iii)} If $|S_0|\ne 2$, then
$\Sigma(P)$ contains no multi-multiple edges.\\
%sub\-dia\-gram $G_2^{(k)}$, for $k>5$.
%
\phantom{{\sc (iii)}} If $|S_0|=2$,
then the edge of $S_0$ has the maximum multiplicity amongst all
edges in $\Sigma(P)$.
%$\Sigma(P)$ contains no  edges
%of multiplicity greater than the multiplicity of the edge of $S_0$.

Suppose that Theorems~A and~B hold for any $d_1$-polytope
satisfying $d_1<d$.
Then there exists a sub\-dia\-gram $S_1\subset\Sigma(P)$
and two vertices $y_0,y_1\in\Sigma(P)$ such that the sub\-dia\-gram
$\l S_0,y_1,y_0,S_1\r $ satisfies the following conditions:

\begin{itemize}
\item[$(1)$]
    $S_0$ and $S_1$ are connected elliptic diagrams,
    $S_0,S_1\ne A_n, D_5$;

\item[$(2)$]
    No vertex of $S_1$ attaches to $S_0$ and %\\
    $|S_0|+|S_1|=d$;

\item[$(3)$]
    $\l y_0,S_1\r $ is either a Lann\'er diagram or one of the four
    diagrams shown in Fig.~\ref{sub_Ess} (in the latter case $y_0$ is the
    marked vertex of the diagram);

\item[$(4)$]
    $\l S_0,y_1\r $ is an indefinite sub\-dia\-gram, and
    one of the following holds:
    \begin{itemize}
    \item[$(4\a)$]
         $y_1$ is not joined to $S_1$,
         and \\
         $\l S_0,y_1\r $
         is either a Lann\'er diagram or one of the four
         diagrams shown in Fig.~\ref{sub_Ess} (in the latter case $y_0$ is the
         marked vertex of the diagram);

    \item[$(4\b)$]
         $y_1$  is a good neighbor of $S_1$,
         and \\ the diagram $\l y_0,S_1\r $ contains no multi-multiple
         edges;

\end{itemize}

\item[$(5)$]
 if $|S_0|\ne 2$, then
$\l S_0,y_1,y_0,S_1\r $  contains no multi-multiple edges;\\
 if $|S_0|=2$, then the edge of $S_0$ has the maximum possible
multiplicity in $\l S_0,y_1,y_0,S_1\r $;

\item[$(6)$]
if $|S_1|=4$ then $S_1$ is a dia\-gram of type $F_4$ or $H_4$;\\
if $|S_1|=3$ then $S_1$ is a dia\-gram of type $B_3$ or $H_3$;\\
 if $|S_1|= 2$, then
 the edge of $S_1$ has the maximum possible multiplicity in
$\l y_0,S_1\r $.

\end{itemize}

\end{lemma}

\begin{figure}[!h]
\begin{center}
\psfrag{m}{\small $m$}
\psfrag{l}{\small $l_1$}
\psfrag{l1}{\scriptsize $$}
\psfrag{lr}{\footnotesize $l_{|S_0|}$}
\psfrag{s}{\small $k$}
\psfrag{y0}{\small $y_0$}
\psfrag{y1}{\small $y_1$}
\psfrag{S0}{\small $S_0$}
\psfrag{S1}{\small $S_1$}
\psfrag{a}{\small $$}
\psfrag{b}{\small $$}
\epsfig{file=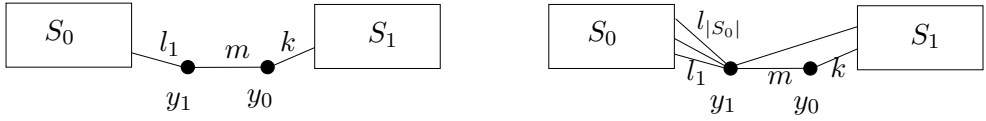,width=0.95\linewidth}
\caption{Diagrams  $\l S_0,y_1,y_0,S_1\r $ satisfying the conditions
$(1)-(6)$: the left one satisfies condition ($4\a$), and the right
one satisfies condition ($4\b$).}
\label{anecdot1}
\end{center}
\end{figure}

Conditions $(1)-(6)$ of the lemma are illustrated in Fig.~\ref{anecdot1}.\\

\vspace{-5pt}

\begin{proof}

We construct the required diagram in several steps.

\begin{itemize}
\item[1.] {\it Analyzing the data.}
Since $S_0$ has no good neighbors,  $\overline S_0=\Sigma_{S_0}$
(see Cor.~\ref{cor_All}).
Denote by $\dim=d-|S_0|$ the dimension of $P(S_0)$.
As a sub\-dia\-gram of $\Sigma(P)$,
the diagram $\Sigma_{S_0}$ contains no dotted edges.
Clearly, $\dim<d$. By the assumption,
 Theorems~A and~B hold for polytopes of dimension less than $d$.
By Prop~\ref{small}, $P(S_0)$ is either a compact simplex
(and $2\le \dim\le 4$) or one
of the Esselmann polytopes (and $\dim=4$).

\item[2.] {\it Choosing $S_1$.}
%and splitting of the ?procedure? into cases (a) and (b).}
We take a sub\-dia\-gram $S_1\subset \overline S_0=\Sigma_{S_0}$ as follows:

If $\dim=4$ then $\overline S_0$ contains a sub\-dia\-gram $S_1$ of type
$F_4$ or $H_4$.

If $\dim=3$ then $\overline S_0$ contains a sub\-dia\-gram $S_1$ of type
$B_3$ or $H_3$.

If $\dim=2$ then $\overline S_0$ contains a sub\-dia\-gram of type
$G_2^{(k)}$, $k\ge 4$, i.e. a multiple edge. We choose $S_1$ as a
diagram $G_2^{(k)}\subset \overline S_0$, where $k$ is maximal in
$\overline S_0$.

Clearly, in all cases conditions $(1)$ and $(2)$ are satisfied.
Notice also, that if $\overline S_0$ contains a multi-multiple edge,
then the diagram $S_1$ has no good neighbors in $\Sigma(P)$.

\item[3.] {\it Choosing $y_0$.}
If $P(S_0)$ is a simplex, then $S_1$ contains all but one vertex
of $\overline S_0$. Let $y_0=\overline S_0\setminus S_1$.

If $P(S_0)$ is an Esselmann polytope,
then it is always possible to choose $y_0\in \overline
S_0\setminus S_1$ such that the diagram $\l  y_0, S_1 \r $
coincides with one of the four diagrams shown in
Fig.~\ref{sub_Ess} (for $y_0$ we take the vertex marked by $y$).

Thus, condition (3) holds.\\

\begin{figure}[!h]
\begin{center}
\psfrag{8}{\scriptsize $8$}
\psfrag{y}{\scriptsize $y$}
\epsfig{file=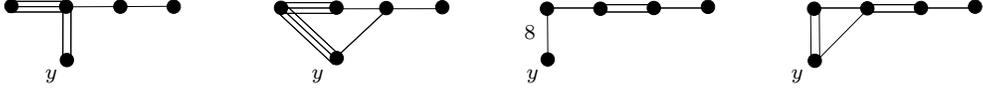,width=0.95\linewidth}
\caption{sub\-dia\-grams of Esselmann diagrams.}
\label{sub_Ess}
\end{center}
\end{figure}

\item[4.] {\it Choosing $y_1$.}
We consider two cases.

\begin{itemize}
\item[($\a$)]
Suppose that $S_1$ has no good neighbors in $\Sigma(P)$.
Then $\overline S_1=\Sigma_{S_1}$
and $P(S_1)$ is either a compact simplex or an Esselmann polytope.
Clearly, $S_0\subset \overline S_1$.
If $P(S_1)$ is a simplex, define $y_1=\overline S_1\setminus S_0$.
If $P(S_1)$ is an Esselmann polytope, we define
$y_1\in \overline S_1$ such that $\l  S_0, y_1\r $
is one of the four diagrams
shown in Fig.~\ref{sub_Ess} (for $y_1$ we take the vertex marked by $y$).
Hence,  $\l  S_0, y_1\r $ satisfies condition $(4\a)$.

\item[($\b$)]
Suppose that $S_1$ has a good neighbor in $\Sigma(P)$.
We choose $y_1$ as one of
good neighbors of $S_1$.
The vertex $y_1$ is connected to $S_1$
by exactly one edge, and this edge is simple.
The vertex $y_1$ might also be connected to any vertex of $S_0$ and
to $y_0$.

Since $S_1$ has a good neighbor and $|S_1|\le d-2$, Cor.~\ref{cor_All}
implies that $S_1\ne F_4$,~$H_4$,~$G_2^{(k)}$ for $k\ge 6$ (in
particular, $S_1$ contains no multi-multiple edge).
Therefore, $\o S_0$ is neither an Esselmann diagram
nor a Lann\'er diagram of order $5$, so $\o S_0$ is a Lann\'er
diagram of order $3$ or $4$.
In the latter case  $\o S_0$ contains no multi-multiple edge.
In the first case, recall that the diagram $S_1$ is chosen as a sub\-dia\-gram
of $\o S_0$ containing the edge of maximal multiplicity.
Since $S_1$ contains no multi-multiple edges,  $\o S_0=\l  y_0, S_1 \r$ does
not contain them, either.
Thus, condition $(4\b)$ holds.

\end{itemize}

\end{itemize}

\noindent
Condition $(5)$ is satisfied by assumption {\sc (iii)} of the lemma,
and condition~$(6)$ is satisfied by the choice of $S_1$, which
completes the proof.

\end{proof}

\begin{lemma}
\label{pr1fin}
The number of diagrams $\l S_0,y_1,y_0,S_1\r $
of signature $(d,1)$, $4\le d \le 8$, satisfying conditions
$(1)-(6)$ of Lemma~\ref{pr1}, is finite.

\end{lemma}

\begin{proof}

Suppose that $S_0\ne G_2^{(k)}$ for $k\ge 6$.
Then by condition $(3)$ the diagram $\l S_0,y_1,y_0,S_1\r $
contains no multi-multiple edges.
Since $|S_0|+|S_1|=d\le 8$, we obtain that
$|\l S_0,y_1,y_0,S_1\r |\le 10$,
and we have finitely many possibilities for the diagram.
%Moreover, the diagram $\l y_0,S_1\r $ is
%either a Lann\'er diagram or the sub\-dia\-grams shown
%in Fig.~\ref{sub_Ess}
%and the sub\-dia\-grams  $\l S_0,y_1\r $ and  $\l y_0,S_1\r $
%are joined by one or two edges only.

Now suppose that $S_0= G_2^{(k)}$, $k\ge 6$.
Since  $|\l S_0,y_1,y_0,S_1\r |=d+2$
and $\sign\l S_0,y_1,y_0,S_1\r =(d,1)$,
we have $\det\l S_0,y_1,y_0,S_1\r =0$.
We consider two cases: either the diagram $S_1$ has a good neighbor
in $\l S_0,y_1,y_0,S_1\r $ or not.

\medskip

\noindent
{\bf Case ($\a$)}: $S_1$ has no good neighbors in $\l S_0,y_1,y_0,S_1\r $.
In this case the sub\-dia\-grams $\l  S_0,y_1 \r $ and
$\l  S_1,y_0 \r $ are either Lann\'er diagrams or diagrams shown in
Fig.~\ref{sub_Ess}. The only edge connecting these diagrams is
$y_0y_1$; we let $m=[y_0,y_1]$ (see Fig.~\ref{anecdot1}).
By Prop.~\ref{loc_product}, we have
$$\det (\l S_0,y_1,y_0,S_1\r,\l y_1,y_0\r) =
\det(\l S_0,y_1\r,y_1) \cdot \det(\l  y_0,S_1 \r,y_0) - \cos^2 \frac{\pi}{m}.$$
Since $\det \l S_0,y_1,y_0,S_1\r = 0$, we obtain
$$\det(\l S_0,y_1\r,y_1) \cdot \det(\l  y_0,S_1 \r,y_0)
=\cos^2 \frac{\pi}{m}.$$
In particular, $$|\det(\l S_0,y_1\r,y_1 )\cdot \det(\l  y_0,S_1 \r,y_0 )|< 1$$
($m=2$ is impossible, since the two indefinite sub\-dia\-grams $\l
S_0,y_1\r $ and $\l  y_0,S_1 \r $ should be joined in
$\Sigma(P)$).
Hence, at least one of
$|\det(\l S_0,y_1\r,y_1 )|$ and $|\det(\l  y_0,S_1 \r,y_0 )|$
is less than 1.

Suppose that $|\det(\l S_0,y_1\r,y_1 )|<1$.
Recall that $S_0=G_2^{(k)}$, $k\ge 6$,
and we have $\det(\l S_0,y_1\r,y_1 )= \d(i,j,k)$,
where $i,j\le k$ by assumption.
Since $|\d(i,j,k)|$ is an increasing function on $i,j,k$,
it is easy to see that if $k\ge 6$, $k\ge i,j$ and $|\d(i,j,k)|<1$,
then $(i,j,k)$ is either $(2,3,7)$ or $(2,3,8)$.
So, $\l S_0,y_1\r $ is either $L_{2,3,7}$ or $L_{2,3,8}$, and
$k\le 8$. Therefore, the diagram
$\l S_0,y_1,y_0,S_1\r $ contains no sub\-dia\-gram $G_2^{(l)}$ for
$l>8$ and we are left with finitely many diagrams.

Suppose that $|\det(\l y_0,S_1\r,y_0 )|<1$.
Since the diagram $\l y_0,S_1\r $ is
either a Lann\'er diagram or one of the diagrams shown
in Fig.~\ref{sub_Ess},
it is easy to check that if $|S_1|>2$ then
$\det(\l y_0,S_1\r,y_0 )>1$.
Therefore, $|S_1|=2$, $d=4$. Again, it is easy to see that there are
only $5$ triples $(i,j,k)$ for which $i,j\le k$ and $|\d(i,j,k)|<1$:
$(i,j,k)=(2,3,7),\; (2,4,5),\; (2,3,8),\; (3,3,4)$ and $(2,5,5)$.
For each of these triples there exist finitely many triples
$(i',j',k')$ satisfying the condition  $|\d(i',j',k')\cdot \d(i,j,k)|<1$.
So, in the case when $S_1$ has no good neighbors in $\l S_0,y_1,y_0,S_1\r $
the lemma is proved.

\medskip

\noindent
{\bf Case ($\b$)}: $y_1$ is a good neighbor of
$S_1$ in $\l S_0,y_1,y_0,S_1\r $.
Note that any edge of  $\l S_0,y_1,y_0,S_1\r $ belongs to either
 $\l S_0,y_1\r $ or  $\l y_1,y_0,S_1\r $.
By Lemma~\ref{loc_sum}, we have
$$ \det(\l S_0,y_1,y_0,S_1\r ,y_1)=
\det(\l S_0,y_1\r ,y_1)+\det(\l y_1,y_0,S_1\r ,y_1)-1.$$
On the other hand,
$$ \det(\l S_0,y_1,y_0,S_1\r ,y_1)=
\frac{\det(\l S_0,y_1,y_0,S_1\r )}{\det(\l S_0,y_0,S_1\r )}=0.$$
Therefore,
$$\det(\l S_0,y_1\r ,y_1)+\det(\l y_1,y_0,S_1\r ,y_1)=1.$$
Since  \ $\l y_0,S_1\r $ \ and \  $\l y_1,y_0,S_1\r $ \ are indefinite
diagrams, \ we obtain that \ $\det(\l y_1,y_0,S_1\r ,y_1)>0$, so
$\det(\l S_0,y_1\r ,y_1)<0$.
Furthermore, $|\l y_1,y_0,S_1\r|=d$, which implies
$$
|\det\l y_1,y_0,S_1\r|  < d!
\eqno (*)
$$
(since the absolute value of each of the summands in the standard
expansion of the determinant does not exceed $1$).
At the same time,
by condition $(4\b)$ the diagram $\l y_0,S_1\r $ contains no
multi-multiple edges,
%sub\-dia\-gram of type $G_2^{(k)}$, $k\ge 6$.
and we have finitely many possibilities for
$\det\l y_0,S_1\r $.
Therefore, there exists a positive constant $M$ such that
$$
M<|\det\l y_0,S_1\r|.
\eqno (**)
$$
Combining $(*)$ and $(**)$, we obtain
$$
0<\det(\l y_1,y_0,S_1\r ,y_1) < \frac {d!}{M},
$$
hence,
$$1-\frac{d!}{M}<\det(\l S_0,y_1\r ,y_1)<0.\eqno (*\!*\!*)$$
Recall that $S_0=G_2^{(k)}$ and that the diagram
$\l S_0,y_1\r $ contains no $G_2^{(l)}$ for $l>k$.
In particular,
$\det(\l S_0,y_1\r,y_1) =\d(i,j,k)$ for some $i,j\le k$. By $(*\!*\!*)$,
we have finitely many possibilities for $k$.
By the assumption,
$\l S_0,y_1,y_0,S_1\r $ contains no sub\-dia\-gram of the type
$G_2^{(l)}$ for $l>k$,
so we have finitely many possibilities for the whole diagram
$\l S_0,y_1,y_0,S_1\r $.

\end{proof}

According to Lemma~\ref{pr1fin},
for each $S_0=G_2^{(k)},B_3,B_4,H_3,H_4,F_4$ we can write down the
complete list $$L(S_0,d)$$ of diagrams  $\l S_0,y_1,y_0,S_1\r $ of
signature $(d,1)$, $4\le d \le 8$, satisfying conditions $(1)-(6)$
of Lemma~\ref{pr1} and containing no parabolic diagrams of order
less than $d$.
Define also a list
\begin{center}
$L(d)=\bigcup\limits_{k=6}^{\infty} L(G_2^{(k)},d).$
\end{center}
By Lemma~\ref{pr1fin}, the list $L(d)$ is also finite.
In view of condition $(4)$ of Lemma~\ref{pr1}, the list
$L(S_0,d)$ naturally splits into two disjoint parts
\begin{center}
$L(S_0,d)=L_{\a}(S_0,d)\cup L_{\b}(S_0,d)$,
\end{center}
where the list $L_{\a}(S_0,d)$ consists of diagrams satisfying
condition $(4\a)$, and the list
$L_{\b}(S_0,d)$ consists of diagrams
satisfying condition $(4\b)$.
Similarly, the list $L(d)$ splits into two parts
\begin{center}
$L_\a(d)=\bigcup\limits_{k=6}^{\infty} L_\a(G_2^{(k)},d)$\quad
and\quad
$L_\b(d)=\bigcup\limits_{k=6}^{\infty} L_\b(G_2^{(k)},d).$
\end{center}

These lists were obtained by a computer.
Usually they are not very short. In what follows we reproduce
some parts of the lists as far as we need.

\vspace{5pt}

\noindent
{\bf Remark.}
It is easy to see that the bounds obtained in the proof of Lemma~\ref{pr1fin}
are not optimal.
In real computations we usually analyze concrete data to reduce calculations.

\vspace{6pt}

%\subsection*{Procedure 2.}

The following lemma is obvious:

\begin{lemma}
For any diagram $\Sigma$ and any constant $C$
the number of diagrams $\l  \Sigma,x\r $
(spanned by $\Sigma$ and a single vertex $x$)
containing no sub\-dia\-grams
$G_2^{(k)}$ for $k>C$ is finite.

\end{lemma}

Hence, for any diagram $\Sigma$, a constant $C$ and dimension $d$,
it is possible to write down a complete list
\begin{center}
$L'(\Sigma,C,d)$
\end{center}
of diagrams   $\l  \Sigma,x\r $ of signature $(d,1)$
containing no sub\-dia\-grams $G_2^{(k)}$ for $k>C$.

Given $\Sigma,C$ and $d$, the list
$L'(\Sigma,C,d)$ can be obtained by a computer.
We reproduce some of these lists as far as we need.
To shorten the computations we use the following:

1) Suppose that $\l  \Sigma, x \r\in L'(\Sigma,C,d)$, and
$|\Sigma|\ge d+1$. Then $|\l  \Sigma, x \r |
\ge d+2$, and $\det\l  \Sigma, x \r =0$.
To check the determinant is faster than to find the signature.
So, first we compute the determinant and in the rare cases when
it vanishes we compute the signature.

2) Suppose that $\Sigma\subset\Sigma(P)$, where $P$ is a simple
hyperbolic $d$-polytope without a pair of disjoint facets.
Suppose that $\Sigma$ contains a connected elliptic sub\-dia\-gram $S\ne A_k,
D_5$. Suppose also that $\overline S\not\subset \Sigma$
(since $|\o S|+|S|>d$,
this always holds if $|\Sigma|\le d+|B|$, where $B$ is the set of bad neighbors
of $S$ in $\Sigma$). In this case there exists $x\in\Sigma(P)\setminus\Sigma$
which is either a good neighbor of $S$ or is not joined to $S$.
Denote by $L'(\Sigma,C,d,S^{(g,n)})$,  $L'(\Sigma,C,d,S^{(g)})$
and  $L'(\Sigma,C,d,S^{(n)})$ the sublists of $L'(\Sigma,C,d)$
which consist of diagrams $\l  \Sigma, x \r $ satisfying the
following conditions $(g,n)$, $(g)$ and $(n)$ respectively:
\begin{itemize}
\item[$(g,n)$]
either $x$ is a {\it good} neighbor of $S$ or $x$ is {\it not}
a neighbor of $S$;

\item[$(g)$]
$x$ is a {\it good} neighbor of $S$;

\item[$(n)$]
$x$ is {\it not} a neighbor of $x$.

\end{itemize}

\noindent
Now we may assume (in the assumptions above) that $\Sigma(P)$ contains a diagram $\l  \Sigma, x \r $
from one of the lists $L'(\Sigma,C,d,S^{(g,n)})$,
$L'(\Sigma,C,d,S^{(g)})$ and  $L'(\Sigma,C,d,S^{(n)})$. This
hugely reduces the computations.

%In the first case, we a left with
%$(const-1)^{n+1-|S|}$ possibilities for the diagram
%$\l  \Sigma, x \r $ (and most of them fail the determinant test).
%If $x$ is a good neighbor of $S$, then there are at most three
%possibility
%to attach $x$ to $S$, and for each of these possibilities we have
%$(const-1)^{n+1-|S|}$ possibilities for the diagram
%$\l  \Sigma, x \r $.

\section{Proof of Theorems~A and~B.}
\label{proof}

The plan of the proof is as follows. We assume that there exists a
simple hyperbolic Coxeter $d$-polytope $P$ with mutually
intersecting facets, and $P$ is not a simplex. Then, using
Theorem~\ref{bor}, Corollary~\ref{ql} and the classification of
Lann\'er diagrams, we find a Coxeter face of $P$ of sufficiently
small codimension.
In view of Theorem~\ref{All}, this face often has no pair of
disjoint facets either.
This enables us to carry out an induction in large dimensions ($d\ge 7$).
In small dimensions (up to $6$) the existence of simplices and
Esselmann polytopes forces us to involve also a computer case-by-case
check based on computations of local determinants.

For a part of the proof (in dimensions $4-6$) we also need
the following lemma.

\begin{lemma}
\label{d+2}
Let $P$ be a simple Coxeter hyperbolic $d$-polytope without a pair of
disjoint facets. If $P$ is neither a simplex nor an Esselmann
polytope nor the polytope shown in Fig.~\ref{t}, then $P$ has at
least $d+3$ facets.

\end{lemma}

\begin{proof}
The lemma follows immediately from the classification of
hyperbolic Coxeter $d$-polytopes with $d+1$ and $d+2$ facets.
The polytopes with $d+1$ facets are simplices,
compact polytopes with $d+2$ facets are either Esselmann polytopes
or simplicial prisms, and the latter have disjoint facets; any
simple non-compact polytope with $d+2$ facets is either a
simplicial prism or the polytope shown in Fig.~\ref{t}.

\end{proof}

\subsection{Dimensions 2 and 3.}
\label{d2,3}

The following lemma does not involve hyperbolic geometry.

\begin{lemma}
\label{2,3}
Let $P$ be a simple $d$-polytope and $d=2$ or $3$.
If $P$ has no pair of disjoint facets then
$P$ is a simplex.

\end{lemma}

\begin{proof}
For $d=2$ the statement is evident.

To prove it for $d=3$ note that any simple 3-polytope different from
simplex has at least one $2$-face which is not a triangle.
Denote such a face by $f$. Let $a$ and $b$ be non-adjacent edges of $f$.
Denote by $f_a$ and $f_b$ the faces of $P$ such that $a=f_a\cap f$
and $b=f_b\cap f$.
By assumption of the lemma $f_a\cap f_b\ne \emptyset$.
Since $P$ is simple, $f_a\cap f_b$ is an edge.
Therefore, the set $\partial P\setminus \{f\cap f_a \cap f_b \}$
has two connected components $M_1$ and $M_2$
(here $\partial P$ is the boundary of $P$).
Each of these components $M_i$ contains
at least one face $m_i$ of $P$,
hence $m_1$ and $m_2$ are two disjoint facets of $P$.

\end{proof}

\subsection{Dimension 4.}
\label{d4}

Lemma~\ref{2,3} does not hold for $4$-polytopes.
Moreover, for any $k\ge 6$ there exists a simple $4$-polytope
with $k$ facets having no pair of disjoint facets. More precisely,
the duals of the cyclic polytopes $C(k,4)$ are simple, have $k$
facets, and any two of its facets intersect in a $2$-face (i.e.
these polytopes are $2$-neighborly); see~\cite{G} for definitions
and details. Furthermore, there are already known seven
Esselmann compact Coxeter hyperbolic $4$-polytopes
with $6$ facets containing no pair of disjoint facets (see
Fig.~\ref{Ess}), and one non-compact $4$-polytope which is
combinatorially equivalent to a product of two simplices (see
Fig.~\ref{t}).

%In this section we suppose that
%$P$ is a Coxeter 4-polytope having no pair of disjoint
%facets.
%Suppose also that $P$ is neither a simplex nor an Esselmann's
%polytope.
%Then by Lemma~\ref{d+2} we immediately obtain $|\Sigma(P)|>6$.
%We consider two possibilities: either  $\Sigma(P)$ contains
%a multi-multiple edge, or it does not (see Prop.~\ref{4_1}
%and~\ref{4_2} respectively).

\begin{prop}
\label{4_1}
Let $P$ be a simple hyperbolic Coxeter $4$-polytope
having no pair of disjoint facets.
If $P$ is not an Esselmann polytope then
 $\Sigma(P)$ contains no multi-multiple edge.

\end{prop}

\begin{proof}
Suppose that $\Sigma(P)$ contains a multi-multiple edge.
Choose $S_0=G_2^{(k)}$, $k\ge 6$, as an edge of
maximal multiplicity in $\Sigma(P)$.
Clearly, $S_0$ has no good neighbors,
so by Lemma~\ref{pr1},
$\Sigma(P)$ contains a sub\-dia\-gram $\l  S_0,y_1,y_0,S_1\r $
from the list $L(4)$.

\begin{figure}[!h]
\begin{center}
\psfrag{q1}{}
\psfrag{q2}{}
\psfrag{6}{\small $6$}
\psfrag{8}{\small $8$}
\psfrag{10}{\small $10$}
\psfrag{a}{\small (a)}
\psfrag{b}{\small (b)}
\psfrag{c}{\small (c)}
\epsfig{file=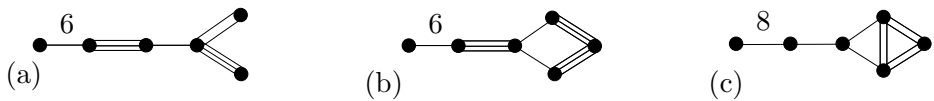,width=0.9\linewidth}
\caption{Intermediate results for $d=4$.}
\label{4new}
\end{center}
\end{figure}

The list  $L_{\a}(4)$ contains two Esselmann diagrams only.
The list $L_{\b}(4)$ contains two Esselmann diagrams and the diagrams shown in
Fig.~\ref{4new}(a),~\ref{4new}(b), and~\ref{4new}(c).
%and those four Esselmann diagrams which contain multi-multiple edges.
By Lemma~\ref{same}, $\Sigma(P)$ contains no Esselmann diagrams.
Hence, we are left with three diagrams shown in
Fig.~\ref{4new}(a),~\ref{4new}(b), and~\ref{4new}(c).
In these cases
$\Sigma(P)$ contains no sub\-dia\-grams $G_2^{(l)}$ for $l>6$,
$6$, and $8$ respectively.
Since none of these diagrams is a diagram of a $4$-dimensional Coxeter
polytope (see Lemma~\ref{d+2}), $\Sigma(P)$ should contain some
diagram from the list $L'(\Sigma,C,4)$, where $\Sigma$ ranges over
the diagrams from
Fig.~\ref{4new}(a),~\ref{4new}(b), and~\ref{4new}(c),
and $C=6$, $6$, and $8$ respectively.
However, these lists are empty: a straightforward computer check
shows that taking any diagram $\l\Sigma,x\r$, where $x$ attaches
to $\Sigma$ by edges of multiplicity at most $C-2$, we obtain
either a superhyperbolic diagram or a diagram with positive
inertia index $\ge 5$ (in fact, we compute the signature only for
those diagrams whose determinant vanishes, see the remark below
the definition of the list $L'(\Sigma,C,d)$).
Thus, we come to a contradiction with the assumption of the
proposition.

\end{proof}

To prove the main result of this section, i.e. Lemma~\ref{4}, we need the following lemma:

\begin{lemma}
\label{3+3, no4}
Let $P$ be a simple hyperbolic Coxeter $4$-polytope having no pair of
disjoint facets. Suppose that $P$ is not a simplex.

(a) Let $v_1,\dots,v_6$ be any six vertices of $\Sigma(P)$.
Then the sub\-dia\-gram spanned by $v_1,\dots,v_6$ contains two
disjoint Lann\'er diagrams of order 3 each.

(b) The order of any Lann\'er sub\-dia\-gram of $\Sigma(P)$ equals 3.

\end{lemma}

\noindent
{\bf Remark.} The lemma involves combinatorics only.
For any simple polytope $P$ we may consider a
``diagram of missing faces'' instead of Coxeter diagram
and missing faces instead of Lann\'er diagrams.

\begin{proof}
Consider a Gale diagram $\G$ of the $4$-polytope $P$.
Denote by $f_1,\dots,f_n$ the facets of $P$.
Then $\G$ is a set of $n$ points at $(n-6)$-dimensional sphere $\S^{n-6}$.
Denote by $b_1,\dots,b_n$ the points corresponding to the facets
$f_1,\dots,f_n$ respectively. Since $P$ is simple, we may assume
that $b_i\ne b_j$ for $i\ne j$.
Denote by $v_i$ the vertex of $\Sigma(P)$ corresponding to a facet $f_i$,
$i=1,\dots,n$.
Consider an $(n-6)$-dimensional plane $\Pi $  spanned
by $b_{7},\dots,b_n$ and the origin.
Again, we may assume that $\Pi$ does not contain points $b_i$ for $i\le 6$.
The hyperplane $\Pi$ separates $\S^{n-6}$ into two hemispheres.
Since $P$ has no disjoint facets, each of the hemispheres
contains at least $3$ points  from $\{b_1,\dots,b_6\}$ (see
Lemma~\ref{missing}).
Hence, three points (say $b_1,b_2,b_3$) belong to one halfspace,
the rest belong to another, which means that $\l v_1,v_2,v_3\r $
and $\l v_4,v_5,v_6\r $ are Lann\'er diagrams (again, see
Lemma~\ref{missing}), and (a) is proved.

To prove (b) suppose that $\l v_k,v_{k+1},v_{k+2},v_{k+3}\r $ is a Lann\'er
diagram. Consider the corresponding points
$b_k,b_{k+1},b_{k+2},b_{k+3}$ in the Gale diagram.
By Lemma~\ref{missing}, there exists an $(n-6)$-plane $\Pi$
through the origin separating these four points. We can rotate the
hyperplane $\Pi$ around the origin until it meets one of the
points $b_1,\dots, b_n$. It cannot meet first any of
$b_k,b_{k+1},b_{k+2},b_{k+3}$ (if $\Pi$ passes through one of
these points then the other three are separated by a plane, so the
four points do not correspond to a Lann\'er diagram).
Hence, $\Pi$ will meet first some point $x_1\in \{b_1,\dots,b_n \}$
distinct from $b_k,b_{k+1},b_{k+2},b_{k+3}$.
Now, we can rotate $\Pi$ around $x_1$ and the origin until $\Pi$ meets
some $x_2\in \{b_1,\dots,b_n \}$,
$x_2\ne b_k,b_{k+1},b_{k+2},b_{k+3}$, and so on.
We have freedom to rotate $\Pi$ until it passes
through $(n-6)$ points $x_1,\dots, x_{n-6}$
(where $x_i\in \{b_1,\dots,b_n \}$,
$x_i\ne b_k,b_{k+1},b_{k+2},b_{k+3}$).
Now $\Pi$ separates $\S^{n-6}$ into two hemispheres:
one contains $4$ points and another contains $n-(n-6)-4=2$ points.
This contradicts the assumption that $P$ have no pair of disjoint
faces. Therefore, no Lann\'er sub\-dia\-gram of $\Sigma(P)$ is of
order $4$. Similarly, it cannot be of order greater than $4$.
Since no Lann\'er sub\-dia\-gram is of order $2$,
we obtain that the order of any Lann\'er diagram equals $3$.

\end{proof}

\begin{lemma}
\label{4}
Let $P$ be a simple hyperbolic Coxeter $4$-polytope.
If $P$ has no pair of disjoint facets, then
$P$ is either a simplex or one of seven Esselmann polytopes
or the polytope shown in Fig.~\ref{t}.

\end{lemma}

\begin{proof}
By Prop.~\ref{4_1}, the diagram
$\Sigma(P)$ contains no multi-multiple edges.
Let $\Sigma \subset \Sigma(P)$ be a sub\-dia\-gram of order $6$
(by Lemma~\ref{d+2}, such a sub\-dia\-gram does exist).
By Lemma~\ref{3+3, no4}, we can assume that $\Sigma=\l S_1, S_2\r $, where
$S_1$ and $S_2$ are Lann\'er diagrams.
There are only $11$ Lann\'er diagrams of order $3$ containing
no edges of multiplicity greater than $5$.
We check all possible pairs of $S_1$ and $S_2$
($66$ possibilities) and connect the vertices of $S_1$ with the
vertices of $S_2$ by edges of all possible multiplicities
($2,3,4,5$ for each of $6$ edges). In all but $39$ cases we obtain
that $\det\l S_1,S_2\r \ne 0$.
Further, $3$ of these $39$ cases correspond to Esselmann diagrams;
one diagram is the diagram of the polytope shown in Fig.~\ref{t};
$4$ diagrams contain parabolic sub\-dia\-grams of order less than $4$;
$11$ of these $39$ diagrams contain Lann\'er sub\-dia\-grams of order $4$, so they
can not be sub\-dia\-grams of $\Sigma(P)$ by Lemma~\ref{3+3, no4}(b).
We are left with $20$ diagrams none of which is a diagram of
Coxeter $4$-polytope: any of them has order $6$, but none of them is
an Esselmann diagram or a diagram of a $4$-prism (see~\cite{Ess}
and~\cite{K}). Therefore, $\Sigma(P)$ contains a sub\-dia\-gram
appearing in one of the lists $L'(\Sigma,5,4)$, where $\Sigma$
ranges over the $20$ diagrams mentioned above.  However, these
lists are empty, and the lemma is proved.

\end{proof}

\subsection{Dimension 5.}
\label{d5}

In this section we suppose that $P$ is a simple hyperbolic
Coxeter $5$-polytope having no pair of disjoint facets.
We also assume that $P$ is not a simplex.

\begin{prop}
\label{5_1}
$\Sigma(P)$ contains neither a sub\-dia\-gram of the type $F_4$
nor a sub\-dia\-gram of the type $H_4$.

\end{prop}

\begin{proof}
Suppose that $\Sigma(P)$ contains a sub\-dia\-gram $S_0=F_4$ or $H_4$.
Then by Cor.~\ref{cor_All}, $\overline S_0=\Sigma_{S_0}$,
and $\Sigma_{S_0}$ contains no dotted edges.
On the other hand, $P(S_0)$ is a $1$-dimensional polytope, i.e. a segment,
so $\Sigma_{S_0}$ should consist of a dotted edge.

\end{proof}

\begin{prop}
\label{5_2}
 $\Sigma(P)$ contains no multi-multiple edges.

\end{prop}

\begin{proof}
Suppose that $\Sigma(P)$  contains a multi-multiple edge.
Choose $S_0=G_2^{(k)}$, $k\ge 6$, as an edge of
maximal multiplicity in $\Sigma(P)$.
Clearly, $S_0$ has no good neighbors, and $\Sigma(P)$ contains a
sub\-dia\-gram appearing in the list $L(5)$.
But all diagrams from the list $L(5)$ contain a sub\-dia\-gram of
the type $H_4$, which contradicts Prop.~\ref{5_1}.

\end{proof}

\begin{prop}
\label{5_3}
$\Sigma(P)$ contains no sub\-dia\-gram of the type $H_3$.

\end{prop}

\begin{proof}
Suppose that $\Sigma(P)$ contains a sub\-dia\-gram $S_0=H_3$.
It follows from Prop~\ref{5_1} that $S_0$ has no good neighbors.
So,  $\Sigma(P)$ contains a
sub\-dia\-gram appearing in the list $L(H_3,5)$.
 The only diagram from the list $L_\a(H_3,5)$ containing
neither a multi-multiple edge nor a sub\-dia\-gram of the type
$H_4$ is shown in Fig.~\ref{5h}(a).
In the list $L_\b(H_3,5)$ there are two diagrams containing
neither a multi-multiple edge nor a sub\-dia\-gram of the type
$H_4$; these diagrams are shown in
Fig.~\ref{5h}(b) and~\ref{5h}(c).

\begin{figure}[!h]
\begin{center}
\psfrag{v}{}
\psfrag{a}{\small (a)}
\psfrag{b}{\small (b)}
\psfrag{c}{\small (c)}
\psfrag{d}{\small (d)}
\psfrag{e}{\small (e)}
\psfrag{f}{\small (f)}
\epsfig{file=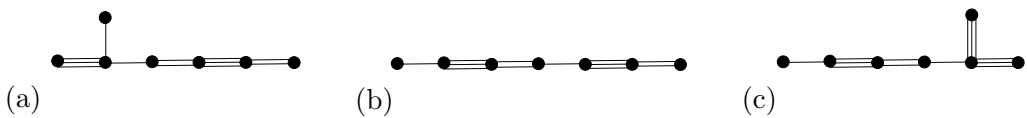,width=0.99\linewidth}
\caption{Intermediate results for $d=5$, $S_0=H_3$.}
\label{5h}
\end{center}
\end{figure}

Consider the diagram $\Sigma$ shown in Fig.~\ref{5h}(a).
By Lemma~\ref{d+2}, $\Sigma$ is not a diagram of a
$5$-polytope. Thus, if $\Sigma(P)$ contains $\Sigma$, then
$\Sigma(P)$ also contains some diagram from the list $L'(\Sigma,5,5)$.
Further, denote by $S$  the sub\-dia\-gram of $\Sigma$ of the type
$B_4$.
Then $\Sigma_S$ is the diagram of a Coxeter $1$-polytope, i.e.
$\Sigma_S$ contains two vertices.
Therefore, $\Sigma(P)$ should contain a diagram from the list
$L'(\Sigma,5,5,S^{(g,n)})$ which happens to be empty.
%
%At the same time, $\Sigma(P)$
%contains no dotted edges.
%It follows that $\Sigma_{S}\ne \o S$
%and by Theorem~\ref{All} at least one of the vertices of
%$\o S$ is a good neighbor of $S$.
%Notice that $\Sigma$ contains only one vertex of $\o S$
%(see Fig.~\ref{5h}(a): any vertex except $v$ either belongs to
%$S$ or is a bad neighbor of $S$).
%Denote by $u$ a good neighbor of $S$. Then $\o S=\l  u,v \r $
%and by Theorem~\ref{All} $[u,v]\ne 2,3$ (otherwise we can not obtain
%a dotted edge in $\Sigma_{S}$). Furthermore, $[u,v]\ne 4$ (otherwise
%$\l  S_2,u,v \r =\wt C_5$ is a parabolic sub\-dia\-gram),
%and $[u,v]\ne 5$ (otherwise $\l  S_2,u,v \r $ contains a
%forbidden diagram $H_4$).
%Therefore, $[u,v]\ge 6$, that is impossible in view of Prop~\ref{5_2}.

Now, consider the diagrams shown in Fig.~\ref{5h}(b) and~\ref{5h}(c).
By Lemma~\ref{d+2}, none of them is a diagram of a
$5$-polytope. Thus, if $\Sigma(P)$ contains one of these two diagrams
(denote it by $\Sigma$), then $\Sigma(P)$ also contains some diagram
from the list $L'(\Sigma,5,5)$.
Furthermore, denote by $S\subset \Sigma$ a diagram of the type
$H_3$ having 2 neighbors in $\Sigma$.
By Prop.~\ref{5_2}, $S$ has no good neighbors in $\Sigma(P)$.
Hence, the diagram $\Sigma(P)$ should contain a diagram
from the list $L'(\Sigma,5,5,S^{(n)})$.
This list turns out to be empty in both cases.
%shown in  Fig.~ref{5h}(b) and (c).

The contradiction shows that the diagrams shown in  Fig.~\ref{5h}(b)
and~\ref{5h}(c) cannot be sub\-dia\-grams of $\Sigma(P)$, which
finishes the proof.

\end{proof}

\begin{prop}
\label{5_4}
$\Sigma(P)$ contains no sub\-dia\-gram of the type $G_2^{(5)}$.

\end{prop}

\begin{proof}
Suppose that $\Sigma(P)\supset S_0=G_2^{(5)}$.
It follows from Prop.~\ref{5_3} that $S_0$ has no good neighbors,
so $\Sigma(P)$ contains a sub\-dia\-gram appearing in the list
$L(S_0,5)$.
However, in the list $L(S_0,5)$ there is no diagram containing
neither a multi-multiple edge nor a sub\-dia\-gram of the types
$H_3$ and $F_4$.

\end{proof}

It follows from Prop.~\ref{5_2} and~\ref{5_4} that
 any multiple edge in $\Sigma(P)$ is a double edge.

\begin{prop}
\label{5_5}
$\Sigma(P)$ contains no sub\-dia\-gram of the type $B_4$.

\end{prop}

\begin{proof}
Suppose that $\Sigma(P)\supset S_0=B_4$.
Then $\o S_0\ne \Sigma_{S_0}$, since $\Sigma_{S_0}$ is a dotted edge
and $\o S_0$ is not.
Let $u$ and $v$ be the vertices of $\o S_0$.
By Theorem~\ref{All}, at least one of
$u$ and $v$ is a good neighbor of $S_0$ (we assume that $u$ is a good
neighbor, so $\l  S_0,u\r$ is either $B_5$, or $\wt B_4$, or
$\wt C_4$).  Suppose that $v$ is not a neighbor of
$S_0$.
Then by Theorem~\ref{All}, $v$ attaches to $u$. If $[u,v]=3$, then
$\l  S_0,u,v\r$ is either of the type $B_6$, or contains a sub\-dia\-gram
of the type $F_4$, or is a quasi-Lann\'er diagram respectively.
If $[u,v]=4$, then $\l  S_0,u,v\r$ is either of the type $\wt C_5$
or contains a sub\-dia\-gram
of the type $\wt C_2$ or, again, is a quasi-Lann\'er diagram.
Now recall that $\Sigma(P)$ contains neither
quasi-Lann\'er diagrams nor elliptic diagrams of order greater
than $5$, and any connected parabolic sub\-dia\-gram of $\Sigma(P)$
should be of order $5$.
Therefore, $v$ is also a good neighbor of $S_0$, the diagram
$\l  S_0,v\r$ is either $B_5$ or $\wt B_4$ or
$\wt C_4$, and the diagram $\Sigma=\l  S_0,u,v\r$ coincides with one of
the diagrams shown in Fig.~\ref{5b}(a)--(g).

\begin{figure}[!h]
\begin{center}
\psfrag{a}{\footnotesize (a)}
\psfrag{b}{\footnotesize  (b)}
\psfrag{c}{\footnotesize  (c)}
\psfrag{d}{\footnotesize  (d)}
\psfrag{e}{\footnotesize (e)}
\psfrag{f}{\footnotesize (f)}
\psfrag{g}{\footnotesize (g)}
\psfrag{h}{\footnotesize (h)}
\psfrag{i}{\footnotesize (i)}
\psfrag{j}{\footnotesize (j)}
\psfrag{k}{\footnotesize (k)}
\epsfig{file=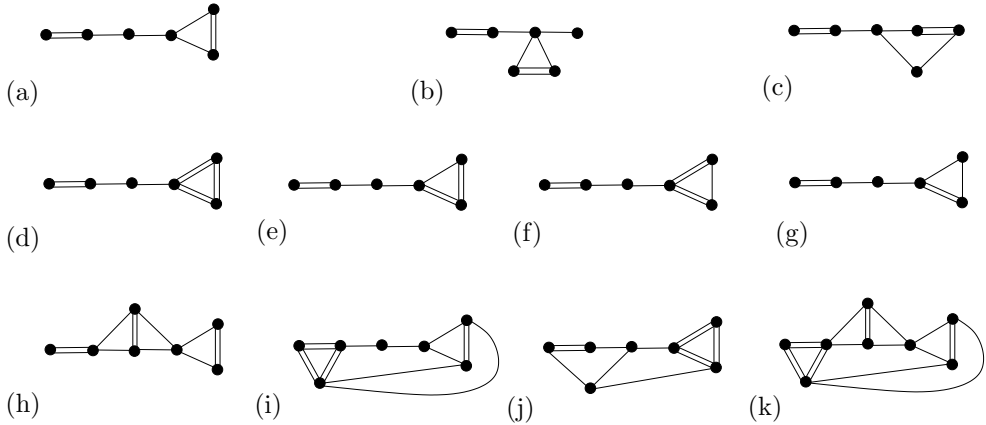,width=0.95\linewidth}
\caption{Intermediate results for $d=5$, $S_0=B_4$.}
\label{5b}
\end{center}
\end{figure}

Since $\Sigma(P)$ contains no multiple edges except for double edges,
$\Sigma(P)$ contains also some diagram from the list
$L'(\Sigma,4,5)$. For the diagram from Fig.~\ref{5b}(a) this list
consists of the two diagrams shown in Fig.~\ref{5b}(h)
and~\ref{5b}(i) (denote these diagrams by $\Sigma_1$ and
$\Sigma_2$). As $\Sigma$ ranges over the diagrams shown in
Fig.~\ref{5b}(b)--(g), the only diagram from a list $L'(\Sigma,4,5)$
which contains neither a sub\-dia\-gram of type $F_4$ nor a
parabolic sub\-dia\-gram of order less than $5$, is that shown in
Fig.~\ref{5b}(j) (denote it by $\Sigma_3$).

Similarly, $\Sigma(P)$ contains some diagram appearing in
either $L'(\Sigma_1,4,5)$ or $L'(\Sigma_2,4,5)$ or
$L'(\Sigma_3,4,5)$.
The latter list is empty, and the former two coincide and
consist of a unique diagram shown in
Fig.~\ref{5b}(k).
The latter diagram contains a sub\-dia\-gram of the type
$ F_4$, and we come to a contradiction.

\end{proof}

\begin{prop}
\label{5_6}
$\Sigma(P)$ contains no sub\-dia\-gram of the type $B_3$.

\end{prop}

\begin{proof}
Suppose that $\Sigma(P)\supset S_0=B_3$.
By Prop.~\ref{5_5}, $S_0$ has no good neighbors and
$\Sigma(P)$ contains some sub\-dia\-gram from the list $L(S_0,5)$.
In the list $L(S_0,5)$ there is no diagram containing neither
sub\-dia\-gram of the type $G_2^{(k)}$, $k\ge 5$, nor sub\-dia\-gram of the
types $B_4$ and $F_4$.

\end{proof}

\begin{prop}
\label{5_7}
$\Sigma(P)$ contains no sub\-dia\-gram of the type $B_2$.

\end{prop}

\begin{proof}
Suppose that $\Sigma(P)\supset S_0=B_2$.
By Prop.~\ref{5_6}, $S_0$ has no good neighbors,
so by Lemma~\ref{2,3},
$\o S_0$ is a Lann\'er diagram of order $4$.
Hence, $\o S_0$ contains a sub\-dia\-gram of the type
either $H_3$ or $B_3$, which is impossible by
Prop.~\ref{5_3} and~\ref{5_6}.

\end{proof}

\begin{lemma}
\label{5}
Let $P$ be a simple hyperbolic Coxeter $5$-polytope.
Then either $P$ has a pair of disjoint facets or $P$
is a non-compact simplex.

\end{lemma}

\begin{proof}
Suppose that $P$ is not a simplex and $P$ has no pair of disjoint facets.
By Prop.~\ref{5_2},~\ref{5_4} and~\ref{5_7},
$\Sigma(P)$ contains no multiple edges. At the same time,
any Lann\'er diagram of order greater than $2$ contains a multiple
edge.
Hence, $\Sigma(P)$ contains no Lann\'er diagram of order greater than $2$.
By Cor.~\ref{ql}, $\Sigma(P)$ contains no quasi-Lann\'er diagram
as well. This means that $\Sigma(P)$ contains a Lann\'er diagram
of order $2$, i.e. a dotted edge.

\end{proof}

\subsection{Dimension 6.}
\label{d6}

In this section we suppose that $P$ is a simple Coxeter
$6$-polytope having no pair of disjoint facets, and $P$
is not a simplex.

\begin{prop}
\label{6_1}
$\Sigma(P)$ contains no multi-multiple edges.

\end{prop}

\begin{proof}
Suppose that $\Sigma(P)$ contains a multi-multiple edge.
Choose $S_0=G_2^{(k)}$, $k\ge 6$, as an edge of
maximal multiplicity in $\Sigma(P)$.
Clearly, $S_0$ has no good neighbors and
$\Sigma(P)$ contains some diagram appearing in the list $L(S_0,6)$.
This list turns out to be empty.

\end{proof}

\begin{prop}
\label{6_2}
$\Sigma(P)$ contains neither sub\-dia\-gram of the type $F_4$ nor
sub\-dia\-gram of the type $H_4$.

\end{prop}

\begin{proof}
Suppose that $\Sigma(P)$ contains a sub\-dia\-gram $S_0$ of the type
either $F_4$ or $H_4$.
Then $S_0$ has no good neighbors and
$\Sigma(P)$ contains a sub\-dia\-gram from the list
$L(S_0,6)$.

The list $L_\a(S_0,6)$ contains a unique diagram $\Sigma$ without
multi-multiple edges. This diagram is shown in Fig.~\ref{6h}(a).
Suppose that $\Sigma\subset\Sigma(P)$.
By Lemma~\ref{d+2}, $\Sigma$ is not a diagram of a $6$-polytope,
hence, $\Sigma(P)$ contains some diagram from the list
$L'(\Sigma,5,6)$.
Further, denote by $S$ the sub\-dia\-gram of $\Sigma$ of the type
$B_5$. Then $\Sigma(P)$ contains also some diagram from the
list $L'(\Sigma,5,6,S^{(g,n)})$. But this list is empty, so
$\Sigma(P)$ contains no sub\-dia\-gram of the type shown in
Fig.~\ref{6h}(a).

\begin{figure}[!h]
\begin{center}
\psfrag{a}{\small (a)}
\psfrag{b}{\small (b)}
\psfrag{c}{\small (c)}
\psfrag{d}{\small (d)}
\psfrag{e}{\small (e)}
\psfrag{f}{\small (f)}
\epsfig{file=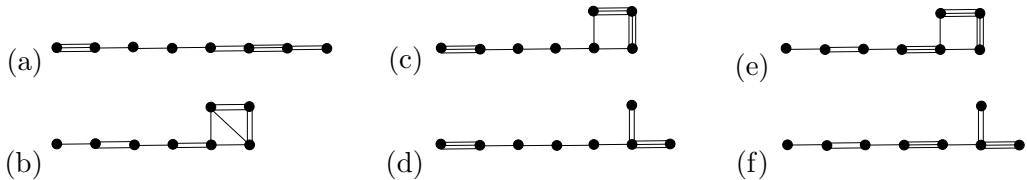,width=0.99\linewidth}
\caption{Intermediate results for $d=6$, $S_0=H_4$ and $F_4$.}
\label{6h}
\end{center}
\end{figure}

The list $L_\b(S_0,6)$ contains five diagrams without
multi-multiple edges. These diagrams are shown in
Fig.~\ref{6h}(b)--(f).
The diagram shown in  Fig.~\ref{6h}(b) contains parabolic
sub\-dia\-grams $\wt C_3$ and $\wt A_2$, which is impossible.
Suppose that $\Sigma(P)$ contains a sub\-dia\-gram $\Sigma$ which is one
of the four diagrams shown in Fig.~\ref{6h}(c)--\ref{6h}(f).
By Lemma~\ref{d+2}, $\Sigma$ is not a diagram of a $6$-polytope,
hence, $\Sigma(P)$ contains some diagram from the list
$L'(\Sigma,5,6)$.
In the cases  Fig.~\ref{6h}(c) and~\ref{6h}(d)
denote by $S$ a sub\-dia\-gram of $\Sigma$ of the type $H_4$
having two neighbors in $\Sigma$.
In the cases  Fig.~\ref{6h}(e) and~\ref{6h}(f) denote by $S$
 a sub\-dia\-gram of $\Sigma$ of the type $H_3$ such that
$S$ is disjoint from the sub\-dia\-gram of the type $F_4$.
Then  $\Sigma(P)$ contains some diagram from the list
$L'(\Sigma,5,6,S^{(g,n)})$.
However, this list is empty in each of the four cases.

\end{proof}

\begin{prop}
\label{6_3}
$\Sigma(P)$ contains no sub\-dia\-gram of the type $H_3$.

\end{prop}

\begin{proof}
Suppose that $\Sigma(P)\supset S_0=H_3$.
By Prop.~\ref{6_2}, $S_0$ has no good neighbors and
$\Sigma(P)$ contains a diagram from the list $L(S_0,6)$.

\begin{figure}[!h]
\begin{center}
\psfrag{a}{\small (a)}
\psfrag{b}{\small (b)}
\psfrag{c}{\small (c)}
\psfrag{d}{\small (d)}
\psfrag{e}{\small (e)}
\psfrag{f}{\small (f)}
\epsfig{file=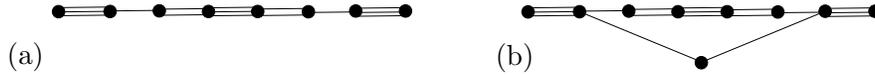,width=0.85\linewidth}
\caption{Intermediate results for $d=6$, $S_0=H_3$.}
\label{6h3}
\end{center}
\end{figure}

In the list  $L_\a(S_0,6)$ there is a unique diagram $\Sigma$ containing
neither multi-multiple edges nor sub\-dia\-gram of the types $H_4$ and $F_4$.
This diagram is shown in Fig.~\ref{6h3}(a).
By Lemma~\ref{d+2}, $\Sigma$ is not a diagram of a $6$-polytope,
so $\Sigma(P)$ contains a sub\-dia\-gram appearing in the list
$L'(\Sigma,5,6)$. Denote by $S$ a sub\-dia\-gram of $\Sigma$ of the type
$B_3$. Then  $\Sigma(P)$ contains a sub\-dia\-gram from the list
$L'(\Sigma,5,6,S^{(g,n)})$.
This list consists of a unique diagram $\Sigma'$ shown in Fig.~\ref{6h3}(b).
The diagram $\Sigma'$ contains a sub\-dia\-gram of the type $H_4$,
which is impossible by Prop~\ref{6_2}.

In the list  $L_\b(S_0,6)$ there is no diagram containing
neither a multi-multiple edge nor a sub\-dia\-gram of the types
$H_4$ and $F_4$.

\end{proof}

\begin{prop}
\label{6_4}
$\Sigma(P)$ contains no sub\-dia\-gram of the type $G_2^{(5)}$.

\end{prop}

\begin{proof}
Suppose that $\Sigma(P)\supset S_0=G_2^{(5)}$.
By Prop.~\ref{6_3}, $S_0$ has no good neighbors, so $\o
S_0=\Sigma_{S_0}$, and $P(S_0)$ is a simple Coxeter 4-polytope
without disjoint facets.
By Lemma~\ref{4}, $\Sigma_{S_0}$ contains either a parabolic
sub\-dia\-gram of the type $\wt C_3$, or a sub\-dia\-gram of the
type $H_4$, or a sub\-dia\-gram of the type $F_4$, which is
impossible by Prop.~\ref{6_2}.

\end{proof}

\noindent
By Prop.~\ref{6_1} and~\ref{6_4}, any multiple edge in
$\Sigma(P)$ is a double edge.

\begin{prop}
\label{6_5}
$\Sigma(P)$ contains no sub\-dia\-gram of the type $B_5$.

\end{prop}

\begin{figure}[!h]
\begin{center}
\psfrag{a}{\small (a)}
\psfrag{b}{\small (b)}
\psfrag{c}{\small (c)}
\psfrag{d}{\small (d)}
\psfrag{e}{\small (e)}
\psfrag{f}{\small (f)}
\psfrag{g}{\small (g)}
\psfrag{h}{\small (h)}
\psfrag{i}{\small (i)}
\epsfig{file=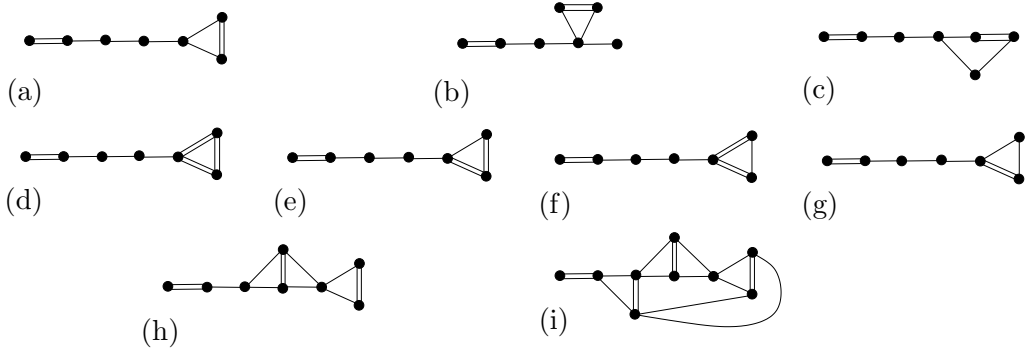,width=0.99\linewidth}
\caption{Intermediate results for $d=6$, $S_0=B_5$.}
\label{6b}
\end{center}
\end{figure}

\begin{proof}
Suppose that $\Sigma(P)\supset S_0=B_5$.
The same argument as in Prop~\ref{5_5} shows that
$\Sigma(P)$ contains a sub\-dia\-gram $\Sigma$ which coincides with one of the
diagrams shown in Fig.~\ref{6b}(a)--(g).
By Lemma~\ref{d+2}, none of these diagrams is a diagram of a $6$-polytope,
so $\Sigma(P)$ contains a diagram from the list $L'(\Sigma,4,6)$.
The union of these lists contains more than $50$ diagrams, but only one of
these diagrams contains neither a sub\-dia\-gram of the type $F_4$
nor a parabolic sub\-dia\-gram of rank less than $5$. This diagram
$\Sigma'$ is shown in Fig.~\ref{6b}(h).
By Lemma~\ref{d+2}, this diagram is not a diagram of a $6$-polytope,
so $\Sigma(P)$ contains a diagram from the list $L'(\Sigma',4,6)$.
The list $L'(\Sigma',4,6)$ consists of a unique diagram $\Sigma''$
shown in Fig.~\ref{6b}(i).
However, the diagram $\Sigma''$ contains a sub\-dia\-gram of the
type $F_4$, which is impossible by Prop~\ref{6_2}.

\end{proof}

\begin{prop}
\label{6_6}
$\Sigma(P)$ contains no sub\-dia\-gram of the type $B_4$.

\end{prop}

\begin{proof}
Suppose that $\Sigma(P)\supset S_0=B_4$.
By Prop.~\ref{6_5}, $S_0$ has no good neighbors and
$\Sigma(P)$ contains a sub\-dia\-gram from the list
$L(S_0,6)$.

In the list  $L_\a(S_0,6)$ there is no diagram containing
neither sub\-dia\-gram $G_2^{(k)}$, $k\ge 5$, nor sub\-dia\-gram
of the types $B_5$ and $F_4$. In the list $L_\b(S_0,6)$ there are
two diagrams containing neither a sub\-dia\-gram $G_2^{(k)}$,
$k\ge 5$, nor a sub\-dia\-gram of the types $B_5$ and $F_4$. These
two diagrams are shown in  Fig.~\ref{6b4}(a) and~\ref{6b4}(b).
Both of these diagrams contain parabolic sub\-dia\-grams of order
$3$, which is impossible.

\begin{figure}[!h]
\begin{center}
\psfrag{a}{\small (a)}
\psfrag{b}{\small (b)}
\psfrag{c}{\small (c)}
\psfrag{d}{\small (d)}
\psfrag{e}{\small (e)}
\psfrag{f}{\small (f)}
\epsfig{file=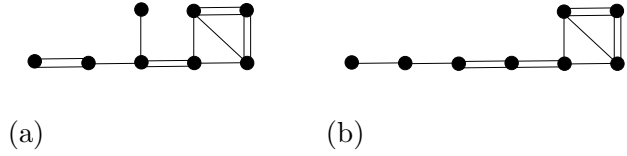,width=0.6\linewidth}
\caption{Intermediate results for $d=6$, $S_0=B_4$.}
\label{6b4}
\end{center}
\end{figure}

\end{proof}

\begin{prop}
\label{6_7}
$\Sigma(P)$ contains no sub\-dia\-gram of the type $B_3$.

\end{prop}

\begin{proof}
Suppose that $\Sigma(P)$ contains a sub\-dia\-gram $S_0=B_3$.
By Prop.~\ref{6_6} $S_0$ has no good neighbors and
$\Sigma(P)$ contains a sub\-dia\-gram from the list $L(S_0,6)$.
However,
in the list $L(S_0,6)$ there is no diagram containing neither a
sub\-dia\-gram $G_2^{(k)}$, $k\ge 5$, nor a sub\-dia\-gram of the
types $B_4$ and $F_4$.

\end{proof}

\begin{prop}
\label{6_8}
$\Sigma(P)$ contains no sub\-dia\-gram of the type $B_2$.

\end{prop}

\begin{proof}
Suppose that $\Sigma(P)\supset S_0=B_2$.
By Prop.~\ref{6_7}, $S_0$ has no good neighbors, and
the proof follows the proof of Prop.~\ref{6_4}.

\end{proof}

\begin{lemma}
\label{6}
Let $P$ be a simple hyperbolic Coxeter $6$-polytope.
Then either $P$ has a pair of disjoint facets or $P$ is a non-compact simplex.

\end{lemma}

\begin{proof}
Suppose that $P$ is not a simplex. By Prop.~\ref{6_1},~\ref{6_4} and~\ref{6_8},
$\Sigma(P)$ contains no multiple edges.
Now the proof follows the proof of Lemma~\ref{5}.

\end{proof}

\subsection{Large dimensions.}
\label{d9}

In this section we assume that $P$ is a simple hyperbolic Coxeter
$d$-polytope  ($d\ge 7$) containing no pair of disjoint facets,
and $P$ is not a simplex. We also assume that $P$
is such a polytope of minimal possible dimension.
We recall that $\Sigma(P)$ contains no quasi-Lann\'er diagrams
(see Cor.~\ref{ql}), so if $S_0\subset\Sigma(P)$ is an elliptic
diagram, $\o S_0=\Sigma_{S_0}$, and $\Sigma_{S_0}$ does not
contain dotted edges, then the dimension of $P(S_0)$ is at most
$4$.

\begin{prop}
\label{7_0}
$\Sigma(P)$ contains no multi-multiple edges.

\end{prop}

\begin{proof}
Suppose that $\Sigma(P)$ contains a sub\-dia\-gram
$S_0=G_2^{(k)}$ for some $k>5$. Then $S_0$ has no good neighbors.
Therefore, $P(S_0)$ is a Coxeter $(d-2)$-polytope without a pair
of disjoint facets, and we contradict our assumptions.

\end{proof}

\begin{prop}
\label{7_1}
$\Sigma(P)$ contains  neither sub\-dia\-gram
of the type $H_4$ nor sub\-dia\-gram of the type $F_4$.

\end{prop}

\begin{proof}
Suppose that $\Sigma(P)$ contains a sub\-dia\-gram $S_0=H_4$ or $F_4$.

For $d=7$ we check the lists $L(S_0,d)$. The union of these lists
for $S_0=H_4$ and $F_4$ consists of four diagrams
$\Sigma_1,\dots,\Sigma_4$ shown in Fig.~\ref{7fh}(a)--(d).
Denote by $S$ a sub\-dia\-gram of $\Sigma_i$ of type $H_4$ having
either two ($i=1,2$) or three ($i=3,4$) bad neighbors.
Since any neighbor of $S$ is bad and none of $\Sigma_i$ is a diagram of
a $7$-polytope, $\Sigma(P)$ contains a sub\-dia\-gram from the list
$L'(\Sigma_i,5,7,S^{(n)})$ for some $i\le 4$. The lists $L'(\Sigma_i,5,7,S^{(n)})$
for $i=1,2,3$ are empty, and the list $L'(\Sigma_4,5,7,S^{(n)})$
consists of a unique diagram $\Sigma_4'$ shown in
Fig.~\ref{7fh}(e).
Again, $\Sigma(P)$ should contain a sub\-dia\-gram from the list
$L'(\Sigma'_1,5,7,S^{(n)})$ for the same $S$. However, this list
is empty.

\begin{center}
\begin{figure}[!h]
\begin{center}
\psfrag{a}{\small (a)}
\psfrag{b}{\small (b)}
\psfrag{c}{\small (c)}
\psfrag{d}{\small (d)}
\psfrag{e}{\small (e)}
\epsfig{file=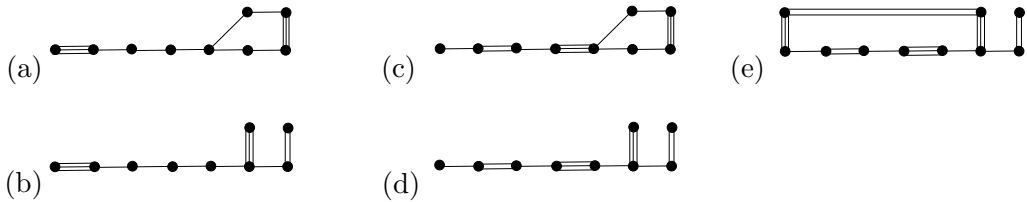,width=0.99\linewidth}
\caption{Intermediate results for $d=7$, $S_0=H_4$ and $F_4$.}
\label{7fh}
\end{center}
\end{figure}
\end{center}

For $d=8$ we check the lists $L(S_0,d)$ which turn
out to be empty.

For $d> 8$ consider the $(d-4)$-polytope $P(S_0)$. By Cor.~\ref{cor_All},
$\Sigma_{S_0}=\overline S_0$. It follows that $\Sigma_{S_0}$ contains
no dotted edges and $P(S_0)$ is a Coxeter $(d-4)$-polytope without
pair of disjoint facets.
If $d=9$ or $d=10$, this contradicts Lemmas~\ref{5} and~\ref{6}
respectively. If $d>10$, this contradicts the assumption that $d$
is the minimal possible dimension of such a polytope.

\end{proof}

\begin{prop}
\label{7_2}
$\Sigma(P)$ contains no sub\-dia\-gram of the type $H_3$.

\end{prop}

\begin{proof}
Suppose that $\Sigma(P)$ contains a sub\-dia\-gram $S_0=H_3$.
By Prop.~\ref{7_1}, $S_0$ has no good neighbors.
Thus,  it follows from Cor.~\ref{cor_All} that
$P(S_0)$ is a Coxeter $(d-3)$-polytope without a pair of disjoint facets.
If $d>7$ then as in Prop.~\ref{7_1} we have a contradiction.

Suppose that $d=7$.
Then $\overline S_0$ is either a Lann\'er diagram of order $5$
or one of the Esselmann diagrams.
In any case, $\overline S_0$ contains either a sub\-dia\-gram of the type
$H_4$ or a sub\-dia\-gram of the type $F_4$, which is impossible
by Prop.~\ref{7_1}.

\end{proof}

\begin{prop}
\label{7_3}
$\Sigma(P)$ contains no sub\-dia\-gram of the type $G_2^{(5)}$.

\end{prop}

\begin{proof}
Suppose that $\Sigma(P)$ contains a sub\-dia\-gram $S_0=G_2^{(5)}$.
Since $\Sigma(P)$ contains no sub\-dia\-gram of the type
$H_3$ (Prop.~\ref{7_2}), $S_0$ has no good neighbors.
Thus,  %it follows from Cor.\ref{Al}(a) that
$P(S_0)$ is a Coxeter $(d-2)$-polytope without a pair of disjoint facets,
and we come to a contradiction.

\end{proof}

\noindent
As a corollary of Prop.~\ref{7_3}, we may assume that all multiple
edges in $\Sigma(P)$ are double edges.

\begin{figure}[!h]
\begin{center}
\psfrag{1}{${\cal L}_1 $}
\psfrag{2}{${\cal L}_2 $}
\psfrag{3}{${\cal L}_3 $}
\psfrag{4}{${\cal L}_4 $}
\psfrag{5}{${\cal L}_5 $}
\epsfig{file=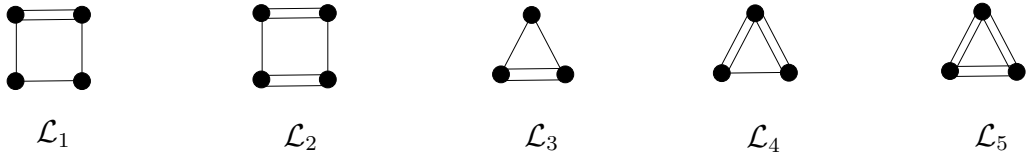,width=0.99\linewidth}
\caption{Notation for some Lann\'er diagrams.}
\label{adm_lan}
\end{center}
\end{figure}

\begin{prop}
\label{7_4}
Any Lann\'er sub\-dia\-gram of $\Sigma(P)$ is one of the
five diagrams shown in Fig.~\ref{adm_lan}.

\end{prop}

\begin{proof}
By the assumption $\Sigma(P)$ contains no Lann\'er diagrams of order $2$.
Thus, the statement follows from the classification of Lann\'er diagrams
and Prop.~\ref{7_3}.

\end{proof}

\begin{prop}
\label{7_5}
If $\Sigma(P)$ contains a sub\-dia\-gram $S=B_3$ or $B_2$ then
$S$ has at least $2$ good neighbors. In addition, for any good
neighbor $u$ of $S$ the diagram $\l S,u\r$ is of the type $B_4$ or
$B_3$, respectively.
\end{prop}

\begin{proof}
Consider the Coxeter polytope $P(S)$. Suppose that $P(S)$ has no pair
of disjoint facets. Then, by assumption, the dimension of $P(S)$
is at most $4$, which means that $d=7$ and $S=B_3$ (see
Lemmas~\ref{d5} and~\ref{d6}). As in the proof of Prop.~\ref{7_2},
$\Sigma_S$ contains a sub\-dia\-gram $\Sigma=H_4$ or $F_4$. By
Cor.~\ref{cor_All}, $\o S$ also contains $\Sigma$, which
contradicts Prop.~\ref{7_1}.

Now we may assume that $P(S)$ has a pair of disjoint facets.
Let $v$ and $u$ be the vertices of $\Sigma_S$
joined by a dotted edge. Denote by $\bar v$ and $\bar u$ the
corresponding vertices of $\Sigma(P)$.
In view of Theorem~\ref{All}, we may assume that one of $v$ and $u$,
say $v$, is a good neighbor of $S$ (otherwise
$\[v,u\]=\[\bar v,\bar u\]\ne \infty$).
Suppose that $u$ is not a neighbor of $S$.
By Prop.~\ref{7_3}, $\[\bar v,\bar u\]\le 4$.
If $\[\bar v,\bar u\]=4$ then $\l S,\bar u,\bar v\r =\wt
C_{4}$ or $\wt C_3$ which are parabolic of small order.
Thus,  $\[\bar v,\bar u\]=2$ or $3$. By item $(2b)$ of Theorem~\ref{All}, we have
$\[v,u\]=2$ or $4$ respectively in contradiction to the assumption that
$\[v,u\]=\infty$. Therefore, $u$ is also a good neighbor of $S$.

By Prop.~\ref{small}, $\l S,u\r$ is not parabolic, and
Prop.~\ref{7_1} implies that $\l S,u\r\ne F_4$, which finishes the
proof.
\smallskip

\end{proof}

\begin{prop}
\label{7_6}
$\Sigma(P)$ contains no sub\-dia\-gram of the type ${\cal L}_1\!$
(see Fig.~\ref{adm_lan}).

\end{prop}

\begin{proof}
Suppose the contrary. Denote the vertices of the sub\-dia\-gram
as shown in Fig.~\ref{prf6}(a).
By Prop~\ref{7_5}, the sub\-dia\-gram $B_3=\l x_1,x_2,x_3\r $ has at least
$2$ good neighbors $y_1$ and $y_2$. By Prop.~\ref{7_5},
$\l x_1,x_2,x_3,y_i\r =B_4$ for $i=1,2$. Clearly, $\[y_1,y_2\]=4$,
otherwise we have either a parabolic sub\-dia\-gram
$\l x_1,x_2,x_3,y_1,y_2\r=\wt B_4$ or a parabolic sub\-dia\-gram
$\l x_3,y_1,y_2\r=\wt A_2$.
Further, $\[x_4,y_i\]\ne 3$ (otherwise $\l x_4,x_3,y_i\r =\wt A_2$),
and  $\[x_4,y_i\]\ne 4$ (otherwise $\l x_2,x_1,x_4,y_i\r =\wt C_3$).
Hence, by Prop.~\ref{7_3} we have $\[x_4,y_i\]=2$ and
$\l x_1,x_2,x_3,x_4,y_1,y_2\r$ is the diagram shown in
Fig.~\ref{prf6}(b).

\begin{figure}[!h]
\begin{center}
\psfrag{x1}{\footnotesize $x_1$}
\psfrag{x2}{\footnotesize $x_2$}
\psfrag{x3}{\footnotesize $x_3$}
\psfrag{x4}{\footnotesize $x_4$}
\psfrag{y1}{\footnotesize $y_1$}
\psfrag{y2}{\footnotesize $y_2$}
\psfrag{z1}{\footnotesize $z_1$}
\psfrag{z2}{\footnotesize $z_2$}
\psfrag{a}{\small (a)}
\psfrag{b}{\small (b)}
\psfrag{c}{\small (c)}
\epsfig{file=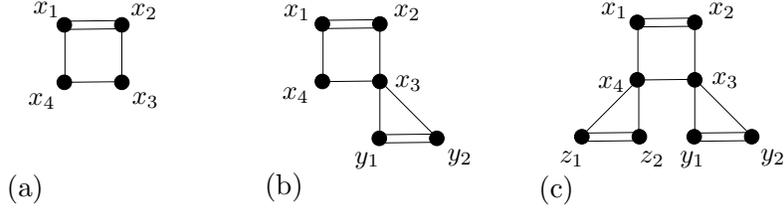,width=0.75\linewidth}
\caption{Notation for the proof of Prop~\ref{7_6}.}
\label{prf6}
\end{center}
\end{figure}

Consider now a pair of good neighbors of the sub\-dia\-gram
$B_3=\l x_4,x_3,x_2\r$ denoting these neighbors by $z_1$ and
$z_2$.
Then $\[z_i,y_j\]\ne 3$ for any $i,j\in \{1,2\}$
(otherwise $\l z_i,x_3,x_4,y_j\r=\wt A_3$).
We also have $\[z_i,y_j\]\ne 4$ for any $i,j\in \{1,2\}$
(otherwise $\l x_1,x_2,x_3,z_i,y_j\r=\wt C_4$ in contradiction
to Prop.~\ref{7_1}).
Thus, $\l x_1,x_2,x_3,x_4,y_1,y_2,z_1,z_2\r$ is the diagram
shown in Fig.~\ref{prf6}(c).
An explicit calculation shows that the sub\-dia\-gram
$\l z_1,z_2,x_4,x_3,y_1,y_2\r$ is superhyperbolic.

\end{proof}

\begin{prop}
\label{7_7}
$\Sigma(P)$ contains no sub\-dia\-gram of the type ${\cal L}_2$.

\end{prop}

\begin{proof}
Suppose the contrary. Denote the vertices of the sub\-dia\-gram
as shown in Fig.~\ref{prf7}(a).
By Prop~\ref{7_5}, the sub\-dia\-gram $B_3=\l x_1,x_2,x_3\r$ has at
least $2$ good neighbors $y_1$ and $y_2$. By Prop~\ref{7_5}, $\l
x_1,x_2,x_3,y_i\r=B_4$.
Clearly, $\[y_1,y_2\]=4$ (see the proof of Prop.~\ref{7_6}).
Further, $\[y_i,x_4]\ne 2$ and $\[y_i,x_4]\ne 4$ (otherwise
we have a parabolic sub\-dia\-gram $\l x_2,x_3,x_4,y_i\r=\wt B_3$ or
$\l y_i,x_4,x_1,x_2\r=\wt C_3$ respectively). Thus, $\[y_i,x_4\]=3$
and $\l x_1,x_2,x_3,x_4, y_1, y_2\r$ is the diagram shown in
Fig.~\ref{prf7}(b).

\begin{figure}[!h]
\begin{center}
\psfrag{x1}{\footnotesize $x_1$}
\psfrag{x2}{\footnotesize $x_2$}
\psfrag{x3}{\footnotesize $x_3$}
\psfrag{x4}{\footnotesize $x_4$}
\psfrag{y1}{\footnotesize $y_1$}
\psfrag{y2}{\footnotesize $y_2$}
\psfrag{z1}{\footnotesize $z_1$}
\psfrag{z2}{\footnotesize $z_2$}
\psfrag{a}{\small (a)}
\psfrag{b}{\small (b)}
\psfrag{c}{\small (c)}
\epsfig{file=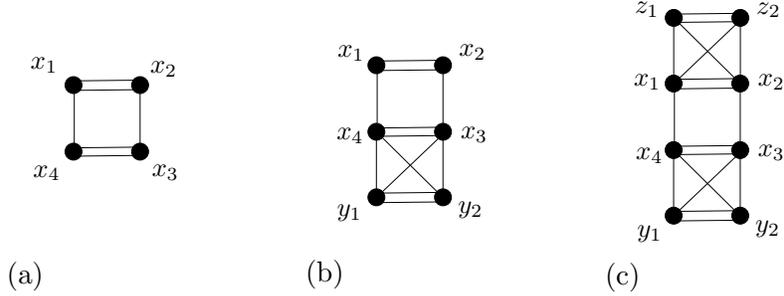,width=0.75\linewidth}
\caption{Notation for the proof of Prop.~\ref{7_7}.}
\label{prf7}
\end{center}
\end{figure}

Consider a pair of good neighbors of the sub\-dia\-gram
$B_3=\l x_4,x_3,x_2\r$ denoting them by $z_1$ and $z_2$.
Then $\[z_i,y_j\]\ne 3$ or $4$ for any $i,j\in \{1,2\}$
(otherwise, we have respectively $\l z_i,x_1,x_4,y_j\r=\wt A_3$ or
${\cal L}_1$).
Thus, $\[z_i,y_j\]=2$ and
$\l x_1,x_2,x_3,x_4,y_1,y_2,z_1,z_2\r$ is the diagram
shown in Fig.~\ref{prf7}(c).
An explicit check shows that the sub\-dia\-gram $\l z_1,z_2,x_2,x_3,y_1,y_2\r$
is superhyperbolic.

\end{proof}

\begin{prop}
\label{7_8}
$\Sigma(P)$ contains no sub\-dia\-gram of the type  ${\cal L}_5$.

\end{prop}

\begin{proof}
Suppose the contrary. Denote the vertices of the sub\-dia\-gram
as shown in Fig.~\ref{prf8}(a).
By Prop.~\ref{7_5}, the sub\-dia\-gram $B_2=\l x_1,x_2\r$
has at least $2$ good neighbors $y_1$ and $y_2$.
We may assume that $\[y_1,x_1\]=2$ and $\[y_1,x_2\]=3$.
Then we have similar conditions for $y_2$:
$\[y_2,x_1\]=2$ and $\[y_2,x_2\]=3$
(otherwise we have $\[y_2,x_1\]=3$ and $\[y_2,x_2\]=2$, so
 the diagram $\l x_1,x_2,y_1,y_2\r$ is either $F_4$ (forbidden
by Prop.~\ref{7_1}) or a cyclic Lann\'er diagram forbidden by
Prop.~\ref{7_6} and~\ref{7_7}). Clearly, $\[y_1,y_2\]=4$, otherwise
either $\l x_1,x_2,y_1,y_2\r=\wt B_3$ or $\l x_2,y_1,y_2\r=\wt A_2$.
Furthermore, $\[y_i,x_3\]\ne 4$ (otherwise $\l y_i,x_3,x_1\r=\wt C_2$),
and $\[y_1,x_3\]=3$ if and only if $\[y_2,x_3\]=3$
(otherwise $\l y_2,y_1,x_3,x_1\r=\wt C_3$). Thus,
$\l x_1,x_2,x_3,y_1,y_2\r$ is one of two diagrams shown in
Fig.~\ref{prf8}(b) and~\ref{prf8}(c).

\begin{figure}[!h]
\begin{center}
\psfrag{x1}{\footnotesize $x_1$}
\psfrag{x2}{\footnotesize $x_2$}
\psfrag{x3}{\footnotesize $x_3$}
\psfrag{x4}{\footnotesize $x_4$}
\psfrag{y1}{\footnotesize $y_1$}
\psfrag{y2}{\footnotesize $y_2$}
\psfrag{z1}{\footnotesize $z_1$}
\psfrag{z2}{\footnotesize $z_2$}
\psfrag{t1}{\footnotesize $t_1$}
\psfrag{t2}{\footnotesize $t_2$}
\psfrag{q1}{\footnotesize $q_1$}
\psfrag{q2}{\footnotesize $q_2$}
\psfrag{a}{\small (a)}
\psfrag{b}{\small (b)}
\psfrag{c}{\small (c)}
\psfrag{d}{\small (d)}
\psfrag{e}{\small (e)}
\psfrag{f}{\small (f)}
\epsfig{file=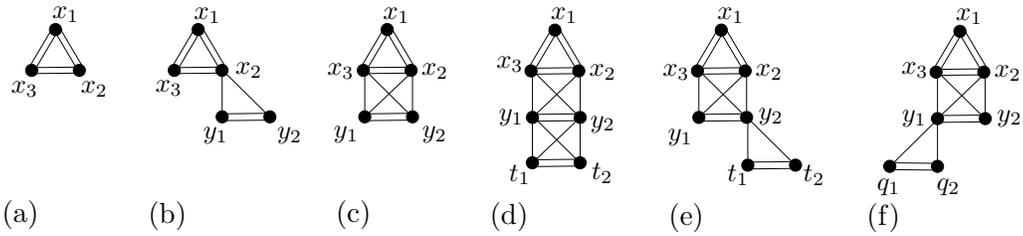,width=0.98\linewidth}
\caption{Notation for the proof of Prop~\ref{7_8}.}
\label{prf8}
\end{center}
\end{figure}

Suppose that $\Sigma(P)$ contains no sub\-dia\-gram of the type shown in
Fig.~\ref{prf8}(c), i.e. $\l x_1,x_2,x_3,y_1,y_2\r$ is the diagram
shown in Fig.~\ref{prf8}(b).
Consider good neighbors $z_1$ and $z_2$ of the sub\-dia\-gram
$B_2=\l x_1,x_3\r$. Without loss of generality, we may assume that
$z_1$ and $z_2$ are neighbors of $x_3$. By the assumption,
$z_1$ and $z_2$ are not neighbors of $x_1$ and $x_2$
(otherwise $\l x_1,x_2,x_3,z_1,z_2\r$ is a sub\-dia\-gram of
the type shown in Fig.~\ref{prf8}(c)).
It follows that the sub\-dia\-gram $\l z_1,x_3,x_2,y_1\r$ is either
$F_4$ (forbidden by Prop.~\ref{7_1}) or a cyclic Lann\'er diagram
(forbidden by Prop.~\ref{7_6} and~\ref{7_7}). We come to a contradiction
which shows that $\Sigma(P)$ contains a sub\-dia\-gram of the type shown in
Fig.~\ref{prf8}(c). We now assume that $\l x_1,x_2,x_3,y_1,y_2\r$
is this sub\-dia\-gram.

Consider two good neighbors $t_1$ and $t_2$ of the sub\-dia\-gram
$B_3=\l x_1,x_2,y_2\r$.
We have

1) $\[t_i,x_3\]=2$ (otherwise either $\l t_i,y_2,x_3\r=\wt A_2$
or $\l x_1,x_3,t_i\r=\wt C_2$);

2) $\[t_1,t_2\]=4$ (otherwise either $\l t_1,t_2,y_2\r=\wt A_2$ or
$\l x_1,x_2,y_2,t_1,t_2\r=\wt B_4$);

3) $\[t_i,y_1\]\ne 4$ (otherwise $\l x_1,x_3,y_1,t_i\r=\wt C_3$);

4) either $\[t_1,y_1\]=\[t_2,y_2\]=3$ or $\[t_1,y_1\]=\[t_2,y_2\]=2$
(otherwise \\ $\l x_1,x_3,y_1, t_1, t_2\r=\wt C_4$).\\
Therefore, $\l x_1,x_2,x_3,y_1,y_2,t_1,t_2\r$ is one of two
diagrams shown in Fig.~\ref{prf8}(d) and~\ref{prf8}(e).
The diagram shown in Fig.~\ref{prf8}(d) is superhyperbolic.
Thus, $\l x_1,x_2,x_3,y_1,y_2,t_1,t_2\r$ is the diagram shown in
Fig.~\ref{prf8}(e).

Consider two good neighbors $q_1$ and $q_2$ of the sub\-dia\-gram
$B_3=\l x_1,x_3,y_1\r$. Reasoning as above shows that the
sub\-dia\-gram $\l x_1,x_2,x_3,y_1,y_2,q_1,q_2\r$ looks like the
diagram shown in Fig.~\ref{prf8}(f).
Then the sub\-dia\-gram $\l q_1,y_1,y_2,t_1\r$ is either $F_4$ (forbidden
by Prop.~\ref{7_1}) or a cyclic Lann\'er diagram (forbidden by
Prop.~\ref{7_6} and~\ref{7_7}).
The contradiction proves the statement.

\end{proof}

\begin{prop}
\label{7_9}
$\Sigma(P)$ contains no sub\-dia\-gram of the type $B_2$.

\end{prop}

\begin{proof}
Suppose the contrary. Let $\l x_1,x_2\r$ be the vertices of $B_2$.
Let $y_1$ and $y_2$ be two good neighbors of $\l x_1,x_2\r$.
Clearly, $\l x_1,x_2,y_1,y_2\r$ is the diagram shown in
Fig.~\ref{prf9}(a).

Let $z_1$ and $z_2$ be two good neighbors of $B_2=\l y_1,y_2\r$.
We have $\[z_i,x_2\]\ne 3$ (otherwise $\l x_2,y_2,z_i\r=\wt A_2$).
If $\[z_1,x_2\]=4$, then $\[z_2,x_2\]=4$ (otherwise
$\l z_1,z_2,x_2\r=\wt C_2$), and $\l z_1,z_2,x_2\r={\cal L}_5$,
which contradicts Prop.~\ref{7_8}.
So, $\[z_i,x_2\]=2$. Furthermore, $\[z_i,x_1\]=2$, otherwise the
cycle $\l x_1,x_2,y_2,z_i\r={\cal L}_1$ or ${\cal L}_2$, which
contradicts Prop.~\ref{7_6} and~\ref{7_7}.
Thus, $\l x_1,x_2,y_1,y_2,z_1,z_2\r$ is the diagram shown in
Fig.~\ref{prf9}(b).

\begin{figure}[!h]
\begin{center}
\psfrag{x1}{\footnotesize $x_1$}
\psfrag{x2}{\footnotesize $x_2$}
\psfrag{x3}{\footnotesize $x_3$}
\psfrag{x4}{\footnotesize $x_4$}
\psfrag{y1}{\footnotesize $y_1$}
\psfrag{y2}{\footnotesize $y_2$}
\psfrag{z1}{\footnotesize $z_1$}
\psfrag{z2}{\footnotesize $z_2$}
\psfrag{t1}{\footnotesize $t_1$}
\psfrag{t2}{\footnotesize $t_2$}
\psfrag{u1}{\footnotesize $q_1$}
\psfrag{u2}{\footnotesize $q_2$}
\psfrag{a}{\small (a)}
\psfrag{b}{\small (b)}
\psfrag{c}{\small (c)}
\psfrag{d}{\small (d)}
\epsfig{file=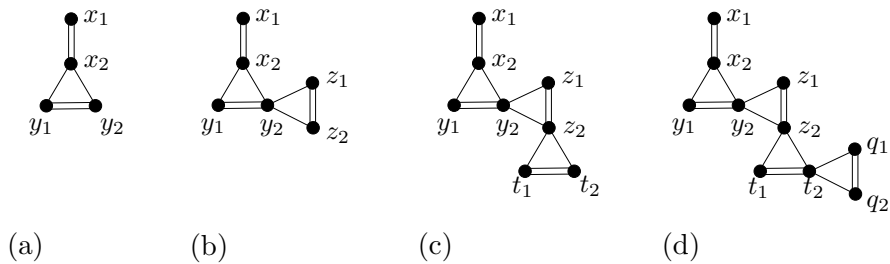,width=0.85\linewidth}
\caption{Notation for the proof of Prop.~\ref{7_9}.}
\label{prf9}
\end{center}
\end{figure}

Let $t_1$ and $t_2$ be two good neighbors of $B_2=\l z_1,z_2\r$.
Repeating the argument above we obtain that
$t_i$ is not connected to $z_1$, $y_2$ and $y_1$.
Moreover, $\[t_i,x_2\]=2$ (see $\l x_2,y_2,z_2,t_i\r$),
and $\[t_i,x_1\]=2$ (otherwise either $\l t_2,x_1,x_2\r=\wt C_2$
or $\l t_i,x_1,x_2,y_1\r=F_4$).
Thus, $\l x_1,x_2,y_1,y_2,z_1,z_2,t_1,t_2\r$ is the diagram shown
in Fig.~\ref{prf9}(c).
This diagram is superhyperbolic, and the proof is complete.

\end{proof}

\noindent
{\bf Remark.} If we consider two good neighbors $q_1,q_2$ of
$B_2=\l t_1,t_2\r$, we obtain a diagram  shown in
Fig.~\ref{prf9}(d), which is evidently superhyperbolic.

\bigskip

%\noindent
Now we are able to finish the proof of the theorems.

\begin{lemma}
\label{9}
Let $P$ be a simple hyperbolic Coxeter $d$-polytope.
%Theorem~\ref{nodots} holds for $d\ge 7$.
If $d>9$ then $P$ has a pair of disjoint facets.
If $6<d\le 9$ then either $P$ has a pair of disjoint
facets or $P$ is a non-compact simplex.
\end{lemma}

\begin{proof}
Suppose that $P$ is not a simplex.
It follows from Prop.~\ref{7_4} and~\ref{7_9} that $\Sigma(P)$ contains
no Lann\'er sub\-dia\-grams of order greater than $2$.
Therefore, it contains a dotted edge, and the lemma is proved.

\end{proof}

\end{document}